%% file: Diplomarbeit.tex
\input{Meta}

\documentclass[12pt,a4paper,twoside]{article}

\usepackage[ngerman, english]{babel}
\usepackage[latin1]{inputenc}
\usepackage[T1]{fontenc}
\usepackage{graphicx}
\usepackage{latexsym}
\usepackage{amsmath,amssymb, amsthm}
\usepackage{ifpdf}
\usepackage[DIV14,BCOR2cm]{typearea}
\usepackage{enumerate}
\usepackage[arrow, matrix, curve]{xy}
\usepackage{ledmac}
\usepackage{makeidx} 
\usepackage{url}
\usepackage{hyperref}
\usepackage[automark]{scrpage2}

\newtheoremstyle{thm}%
        {10pt}
        {10pt}
        {\itshape}
        {}
        {\bfseries}
        {}
        {2mm}
        {}
\newtheoremstyle{def}%
        {10pt}
        {10pt}
        {\normalfont}
        {}
        {\bfseries}
        {}
        {2mm}
        {}
\theoremstyle{def}
\newtheorem{Definition}{Definition}[section]
\newtheorem{Example}[Definition]{Example}
\newtheorem{Remarkc}[Definition]{Remark}
\newtheorem*{Remark}{Remark}
\newtheorem*{Proof}{Proof}
\theoremstyle{thm}
\newtheorem{Theorem}[Definition]{Theorem}
\newtheorem{Lemma}[Definition]{Lemma}
\newtheorem{DefinitionLemma}[Definition]{Definition-Lemma}
\newtheorem{Corollary}[Definition]{Corollary}
\numberwithin{equation}{section}

\makeindex

\begin{document}
\include{Deckblatt}
\newpage
\thispagestyle{empty}
\mbox{}
\newpage
\setcounter{page}{1}
\tableofcontents
\newpage
\mbox{}
\newpage
\section{Introduction}
A global action is an algebraic analogue of a topological space. It consists of group actions $G_\alpha\curvearrowright X_\alpha,~ (\alpha\in \Phi)$, which fulfill a certain compatibility condition (see Definition 2.1). In Section 2, immediately below we define three kinds of morphisms between global actions namely a general notion of morphism, then regular morphisms and finally normal morphisms. It turns out that every regular morphism is normal. Further we define the product of two global actions and show how one can put a global action structure on the set $Mor(X,Y)$ of all morphisms from a global action X to a global action Y. The main theorem of Section 2 is the Exponential Law 2.31, which ensures that one can develop a good homotopy theory for global actions. In Section 3 we define the notion of path in a global action X, introduce the notion of homotopy of morphisms and define the fundamental group $\pi_1(X_*)$ of a global action X at a fixed point $*$. Section 4 introduces the notion of covering morphism. We show important lifting properties and as the main result of this section Theorem 4.8, which establishes a Galois type correspondence between connected coverings of a given connected global action and subgroups of the fundamental group of that action.\\
 
\section{Basic Definitions}
Recall that if X denotes a set and G a group then a \textit  {(left) group action} or just \textit {(left) action}
of G on X is a rule which associates to each pair $(g,x)$ of elements $g\in G$ and $x\in X$ an element $gx\in X$ such that if $g_1, g_2\in G$ then $g_2(g_1x)=(g_2g_1)x$ and such that if 1 denotes the identity element in G then $1x=x$. An action of G on X is often denoted by $G\curvearrowright X$.\ A \textit {morphism} $G\curvearrowright X \rightarrow G'\curvearrowright X'$ of left actions is a pair $(h,s)$ consisting of a group homomorphism $h:G\rightarrow G'$ and a function $s:X \rightarrow X'$ such that for any $g\in G$ and $x\in X$ $s(gx)=h(g)s(x)$.

The notion of a right action of G on X and a morphism of right actions is defined similarly .
\Definition[global action]{
Let X denote a set. A (left) \textit {global action}\index{global action} on X consists of a set $\{G_\alpha \curvearrowright X_\alpha| \alpha\in\Phi\}$ of group actions $G_\alpha\curvearrowright X_\alpha$ where each $X_\alpha\subseteq X$ and $\Phi$ is an indexing set equipped with a reflexive relation $\leq$ such that the following holds:\
\begin {enumerate}[(a)]
\item If $\alpha\leq\beta$ then $X_\alpha\cap X_\beta$ is stable under the action of $G_\alpha$, i.e. if $g\in G_\alpha$ and $x\in X_\alpha\cap X_\beta$ then $gx\in X_\alpha\cap X_\beta$. Thus $G_\alpha \curvearrowright(X_\alpha\cap X_\beta)$.
\item Suppose $\alpha\leq\beta\in \Phi$. Then as part of the structure of a global action there is a group homomorphism $G_{\alpha\leq\beta}:G_\alpha\rightarrow G_\beta$ such that the pair $(G_{\alpha\leq\beta}, i_{X_\alpha\cap X_\beta, X_\beta})$ is a morphism $G_\alpha\curvearrowright(X_\alpha\cap X_\beta)\rightarrow G_\beta\curvearrowright X_\beta$ of group actions where $i_{X_\alpha\cap X_\beta, X_\beta}:X_\alpha\cap X_\beta\rightarrowtail X_\beta$ denotes the canonical inclusion. Such a homomorphism $G_{\alpha\leq\beta}$ is called a \textit{structure homomorphism}. It is assumed that $G_{\alpha\leq\alpha}:G_\alpha\rightarrow G_\alpha$ is the identity map.
\end{enumerate}
The set X equipped with a global action is called a \textit {global action space} or simply a \textit {global action}. The sets $X_\alpha$ are called \textit{local sets}, the groups $G_\alpha$ \textit{local groups}, and the actions $G_\alpha\curvearrowright X_\alpha$ \textit{local actions}.
}
\Remark{
\begin {enumerate}[(a)]
\item In this thesis it will additionally be assumed that $X=\bigcup\limits_{\alpha\in\Phi}X_\alpha$ .
\item Since it is not assumed that the relation on $\Phi$ is transitive, it follows that if $\alpha\leq\beta$ and $\beta\leq\gamma$ then it does not follow automatically that $\alpha\leq\gamma$, but it is also not forbidden that $\alpha\leq\gamma$. Suppose $\alpha\leq\gamma$. Then $G_{\alpha\leq\beta}$, $G_{\beta\leq\gamma}$ and $G_{\alpha\leq\gamma}$ are structure homomorphisms but it is not assumed that $G_{\beta\leq\gamma}G_{\alpha\leq\beta}=G_{\alpha\leq\gamma}$, although this equality is also not forbidden. However, in many situations which arise naturally, the relation on $\Phi$ is transitive and the composition rule $G_{\beta\leq\gamma}G_{\alpha\leq\beta}=G_{\alpha\leq\gamma}$ holds.
\end {enumerate}
}\hfill

Below are two important examples of global actions.
\Example{
Let G denote a group and let $\Phi$ denote a set which is indexing a set $\{G_\alpha|\alpha\in\Phi\}$ of subgroups $G_\alpha$ of G, i.e. if $G_\alpha=G_{\alpha'}$ then $\alpha=\alpha'$, which is closed under taking intersections, i.e. if $\alpha,\beta\in\Phi$, $\exists$ a $\gamma\in\Phi\backepsilon G_\alpha\cap G_\beta=G_\gamma$. Equip $\Phi$ with the reflexive, transitive relation defined by $\alpha\leq\beta\Leftrightarrow G_\alpha\subseteq G_\beta$. If $\alpha\leq\beta$ let $G_{\alpha\leq\beta}$ denote the canonical inclusion homomorphism $G_\alpha\rightarrowtail G_\beta$. Let $X=G$ and for each $\alpha\in\Phi$ let $X_\alpha=X$. Let $G_\alpha$ act on $X_\alpha$(=X) by left (or right) multiplication. Then $\{G_\alpha\curvearrowright X_\alpha |\alpha\in\Phi\}$ is a global action on X. This example is called a \textit{standard single domain} global action.
}
\Example[line action\index{line action}]{
Let $X=\mathbb{Z}$, $\Phi=\mathbb{Z}\cup\{*\}$ and $X_n=\{n, n+1\}$ if $n\in\mathbb{Z}$ and $X_*=\mathbb {Z}$. Let $G_n\cong\mathbb{Z}/2\mathbb{Z}$ if $n\in\mathbb{Z}$, $G_*=1$ and let $G_n\curvearrowright \{n, n+1\}$ be the group action such that the non-trivial element of $G_{\{n,n+1\}}$ exchanges the elements n and n+1. Let the only relations in $\Phi$ be $*\leq n$ for all $n\in\mathbb{Z}$. Let $G_{*\leq n}:\ \{1\}\rightarrow G_{\{n,n+1\}}$ denote the unique group homomorphism. This example is called the \textit{line action}.
}\\

%

We prepare now for the notion of a morphism of global actions.
\Definition[local frame\index{local frame}]{
Let $\{G_\alpha\curvearrowright X_\alpha|\alpha\in\Phi\}$ denote a global action on the set X and let $\alpha\in\Phi$. A \textit{local frame} at $\alpha$ or simply an $\alpha$-\textit{frame} is a finite set $\{x_0,...,x_p\}\subseteq X_\alpha$ such that $G_\alpha$ acts transitively on $\{x_0,..., x_p\}$, i.e. for each $i\ (1\leq i\leq p)$, there is an element $g_i\in G_\alpha\backepsilon g_i x_0=x_i$. 
}
\Definition[morphism \index{morphism} of global actions]{
Let $\{G_\alpha\curvearrowright X_\alpha|\alpha\in\Phi\}$ and $\{H_\beta\curvearrowright Y_\beta|\beta\in\Psi\}$ be global actions on X and Y, respectively. Then a \textit{morphism} $f:X\rightarrow Y$ of global actions is a function f which preserves local frames, i.e. if $\{x_0,...,x_p\}$ is a local frame in X then $\{fx_0, ..., fx_p\}$ is a local frame in Y.
}

\Remark{If $f:X\rightarrow Y$, $g:Y\rightarrow Z$ are morphism of global actions then obviously their composition $gf:Y\rightarrow Z$ is a morphism of global actions.}

\Example {Let A be a gobal action on a space X, $B=\{H_\beta\curvearrowright Y_\beta|\beta\in\Psi\}$ be a global action on a space Y such that $\exists \beta \in\Psi$ and an element $y\in Y_\beta$ such that $H_\beta y=Y_\beta$, i.e all elements of $Y_\beta$ are in the orbit of some element $y\in Y_\beta$. Then any function $f:X\rightarrow Y$ whose image lies in $Y_\beta$ is a morphism of global actions.
}
\Example {This example is a special case of Example 2.6. Define the global action B as follows: There is only one group $G_1=\mathbb{Z}$, which acts on itself by left addition (of course $G_{1\leq1}$ is the identity map $id_\mathbb{Z}:\mathbb{Z}\rightarrow\mathbb{Z}$). Then the function $id_\mathbb{Z}$ is a morphism from the global action space $\mathbb{Z}$ equipped with the line action defined in Example 2.3 to $\mathbb{Z}$ equipped with the global action B.
}

\Definition[regular morphism \index{morphism!regular -}]{Let $\{G_\alpha\curvearrowright X_\alpha|\alpha\in\Phi\}$ denote a global action on X and $\{H_\beta\curvearrowright Y_\beta|\beta\in\Psi\}$ a global action on Y. A \textit{regular} morphism $\eta:X\rightarrow Y$ is a triple $(\iota, \kappa, \lambda)$ where 
\begin {enumerate}[(a)]
\item $\iota:\Phi\rightarrow\Psi$ is a function such that $\iota(\alpha)\leq\iota(\beta)$ whenever $\alpha\leq\beta$, i.e. $\iota$ is relation preserving.
\item $\kappa$ is a rule which assigns to each $\alpha \in\Phi$ a group homomorphism $\kappa_\alpha:G_\alpha\rightarrow H_{\iota(\alpha)}$ such that if $\alpha\leq\beta$ then the diagram \[\begin{xy}\xymatrix{G_\alpha \ar[r]^{\kappa_\alpha} \ar[d]_{G_{\alpha\leq\beta}}&H_{\iota(\alpha)} \ar[d]^{H_{\iota(\alpha)\leq\iota(\beta)}} \\G_\beta \ar[r]_{\kappa_\beta} &H_{\iota(b)}}\end{xy}\] commutes. (Thus, if the relations on $\Phi$ and $\Psi$ are transitive (and therefore $\Phi$ and $\Psi$ are categories in the obvious way) and the composition rules $G_{\beta\leq\gamma}G_{\alpha\leq\beta}=G_{\alpha\leq\gamma}$ and $H_{\beta'\leq\gamma'}H_{\alpha'\leq\beta'}=H_{\alpha'\leq\gamma'}$ hold then G and H are functors on respectively $\Phi$ and $\Psi$, with values in groups, $\iota$ is a functor $\Phi\rightarrow\Psi$, and $\kappa$ defines a natural transformation $G_{(\ )}\rightarrow H_{(\ )}$ of functors defined on $\Phi$.) The rule $\kappa$ will be called a \textit{natural transformation} $\kappa:G_{(\ )}\rightarrow H_{(\ )}$ from $G_{(\ )}$ to $H_{(\ )}$ even though $\Phi$ and $\Psi$ might not even be categories. 
\item $\lambda:X\rightarrow Y$ is a function such that $\lambda(X_\alpha)\subseteq Y_{\iota(\alpha)}\ \forall\alpha\in\Phi$
\item For all $\alpha\in\Phi$ the pair $(\kappa_\alpha, \lambda|_{X_\alpha}):G_\alpha\curvearrowright X_\alpha\rightarrow H_{\iota(\alpha)} \curvearrowright Y_{\iota(\alpha)}$ is a morphism of group actions.
\end{enumerate}
A regular morphism $\eta=(\iota,\kappa,\lambda):X\rightarrow Y$ is called an \textit{extension} of a morphism $f:X\rightarrow Y$ if $\lambda=f$. A regular morphism $\eta=(\iota,\kappa,\lambda):X\rightarrow Y$ is called a \textit{regular isomorphism} if there is a regular morphism $\eta'=(\iota',\kappa',\lambda'):Y\rightarrow X$ such that $\iota'$ is inverse to $\iota$, $\lambda'$ is inverse to $\lambda$ and $\kappa'_{\iota(\alpha)}$ is inverse to $\kappa_\alpha$. Such a regular morphism $\eta'$ is called the \textit{inverse} of $\eta$. A regular morphism $\eta=(\iota,\kappa,\lambda):X\rightarrow Y$ is called a \textit{weak regular isomorphism} if there is a regular morphism $\eta'=(\iota',\kappa',\lambda'):Y\rightarrow X$ such that $\lambda'$ is inverse to $\lambda$.
} 

\Remark{Condition 2.8(d) implies that a regular morphism is a morphism of global actions (to be more precise if $(\iota, \kappa, \lambda)$ is a regular morphism then $\lambda$ is a morphism of global actions).
\\

To develop a homotopy theory for morphisms of global actions it is necessary to put the structure of a global action on the set $Mor(X,Y)$ of all morphisms from a global action X to a global action Y.
\Definition {Let $\{G_\alpha\curvearrowright X_\alpha|\alpha\in\Phi\}$ denote a global action on X and $\{H_\beta\curvearrowright Y_\beta|\beta\in\Psi\}$ a global action on Y. Define a global action on $Mor(X,Y)$ \index{global action!morphism space of two -s} as follows. The index set $\Theta$ for the global action to be put on Mor(X,Y) is $\Theta=\{\beta:X\rightarrow\Psi|\beta\ \textnormal{an arbitrary function}\}$. The relation on $\Theta$ is defined by $\beta\leq\beta'\Leftrightarrow\beta(x)\leq\beta'(x)\forall x\in X$. If $\beta\in\Theta$ then the local set $Mor(X,Y)_\beta$ is defined by $Mor(X,Y)_\beta=\{f\in Mor(X,Y)|$ 
\begin {enumerate}[(a)]
\item $f(x)\in Y_{\beta(x)}\ \forall x\in X$
\item if $\{x_0,...x_p\}\subseteq X$ is a local frame in X then $\{f(x_0),...,f(x_p)\}$ is a local b-frame in Y for some $b\in\Psi$ such that $b\geq\beta(x_i)\ (0\leq i\leq p)\}$.
\end{enumerate}
If $\beta\in\Theta$ then the local group $J_\beta$ at $\beta$ is defined by $J_\beta=\prod \limits_{x\in X}H_{\beta(x)}$ and the action of $J_\beta$ on $Mor(X,Y)_\beta$ is given by \begin{align*}J_\beta\times Mor(X,Y)_\beta&\rightarrow Mor(X,Y)_\beta\\ (\sigma,  f)&\mapsto \sigma f\end{align*} where \begin{align*}\sigma f:X&\rightarrow Y\\ x&\mapsto \sigma_x (f(x))\end{align*} and \begin{align*}\sigma_x \textnormal{ is the x-coordinate of } \sigma \in \prod \limits_{x\in X}H_\beta(x).\end{align*} One has to check that the function $\sigma f$ defined above is an element of $Mor(X,Y)_\beta$, but the proof is straightforward. If $\beta,\beta'\in\Theta$ and $\beta\leq\beta'$ then the structure morphism $J_{\beta\leq\beta'}$ is defined by $J_{\beta\leq\beta'}=\prod \limits_{x\in X} H_{\beta(x)\leq\beta'(x)}:J_\beta=\prod \limits_{x\in X} H_{\beta(x)}\rightarrow J_{\beta'}=\prod \limits_{x\in X} H_{\beta'(x)}$.

One checks straightforward, but with a little effort that $(\Theta, J_{(-)},$ $ Mor(X,Y)_{(-)})$ defines a global action structure on the set Mor(X,Y).
}\\

The following lemma is easy to check. 
\Lemma{Let X and Y be global actions. If $g:Z\rightarrow X$ is a morphism of global actions then \begin{align*}Mor(g, 1_Y):Mor(X,Y)&\rightarrow Mor(Z,Y)\\ f&\mapsto fg\end{align*}is a morphism of global actions.}\\

\normalfont Unfortunately it is not necessarily the case that if Z and X are global actions and $g:X\rightarrow Y$ a morphism of global actions then $Mor(1_z,g)$ defines a morphism \begin{align*}Mor(Z,X)&\rightarrow Mor(Z,Y)\\ f&\mapsto gf\end{align*} of global actions. If $Mor(1_z,g)$ is a morphism of global actions then we say that g is Z-normal. We make this now a definition.
\Definition[Z-normal morphism \index{morphism!Z-normal -}]{Let Z be a global action. A morphism $g:X\rightarrow Y$ of global actions is called Z-\textit{normal}, if for any local frame $\{f_0, f_1, ..., f_p\}$ at $\beta$ where $\beta\in\Theta_{(Z,X)}$ there is a $\gamma\in\Theta_{(Z,Y)}$ such that $\{gf_0, gf_1, ..., gf_p\}$ is a local frame at $\gamma$. 
}
\Lemma{
Let Z be a global action and $g:X\rightarrow Y$ be a morphism of global actions. Then the function \begin{align*} Mor(1_Z,g):Mor(Z,X)&\rightarrow Mor(Z,Y)\\ f&\mapsto gf\end{align*} is a morphism of global actions if and only if g is Z-normal.
}
\Proof{$Mor(1_Z,g)$ is a morphism of global actions $\Leftrightarrow$ for all local frames $\{f_0, f_1, ..., f_p\}$ in $Mor(Z,X)$ $\{Mor(1_Z,g)f_0, Mor(1_Z,g)f_1, ..., Mor(1_Z,g)f_p\}=\{gf_0, gf_1, ..., gf_p\}$ is a local frame in $Mor(Z,Y)$ $\Leftrightarrow$ g is Z-normal.\hfill$\Box$
}
\Definition[normal morphism\index{morphism!normal -}, $\infty$-normal morphism\index{morphism!$\infty$-normal -}]{Let $g:X\rightarrow Y$ be a morphism of global actions. g is called \textit{normal}, if it is Z-normal for all global actions Z. g is called \textit{$\infty$-normal}, if for any finite set of global actions $Z_1,...Z_n$ the map $Mor(1_{Z_n},Mor(1_{Z_{n-1}},...,Mor(1_{Z_1},g))...):Mor(Z_n,Mor(Z_{n-1},...,Mor(Z_1,X))...)\rightarrow Mor(Z_n,$ $ Mor(Z_{n-1}, ...,Mor(Z_1,Y))...)$ is a morphism of global actions.}\\

The following lemma is obvious and useful.
\Lemma{If $g:X\rightarrow Y$ is an $\infty$-normal morphism of global actions then for any global action Z $Mor(1_Z,g):Mor(Z,X)\rightarrow Mor(Z,Y)$ is $\infty$-normal.}
\Definition{A morphism $g:X\rightarrow Y$ of global actions is called a \textit{Z-normal} (resp. \textit{normal}, \textit{$\infty$-normal}) \textit{isomorphism}, if the map g is bijective and $g^{-1}$ is Z-normal(resp. normal, $\infty$-normal).
}
\Definition{A global action X is called \textit{Z-normal}\index{global action!Z-normal -} (resp. \textit{normal}\index{global action!normal -} , \textit{$\infty$-normal}\index{global action!$\infty$-normal -} ), if every morphism whose domain is X is Z-normal (resp. normal, $\infty$-normal). A global action Y is called \textit{Z-conormal} (resp. \textit{conormal}, \textit{$\infty$-conormal}), if every morphism whose codomain is Y is Z-normal (resp. normal, $\infty$-normal).
}

\Lemma{
If $(\iota,\kappa,\lambda):X\rightarrow Y$ is a regular morphism then for any global action Z, the function $Mor(1_Z,\lambda):Mor(Z,X)\rightarrow Mor(Z,Y)$ extends to a regular morphism. 
}
\Proof{We shall define a relation preserving map $\iota':\Theta_{(Z,X)}\rightarrow \Theta_{(Z,Y)}$ and a natural transformation $\kappa':J_{(Z,X))}\rightarrow J_{(Z,Y)}$ and show that the triple $(\iota',\kappa', Mor(1_Z,g))$ is a regular morphism $Mor(Z,X)\rightarrow Mor(Z,Y)$ of global actions. Define $\iota'$ by $\iota'(\beta)=\iota\beta\ \forall\beta\in\Theta_{(Z,X)}$ and $\kappa'_\beta=\prod\limits_{z\in Z}\kappa_{\beta(z)}\ \forall\beta\in\Theta_{(Z,X)}$. We show now that $(\iota',\kappa', Mor(1_Z,g))$ is a regular morphism. 
\begin {enumerate}[(a)]
\item Let $\beta,\gamma\in\Theta_{(X,Y)}, \beta\leq\gamma$ i.e. $\beta(z)\leq\gamma(z)$ for all z in Z. Then $(\iota'(\beta))(z)=(\iota\beta)(z)=\iota(\beta(z))\leq\iota(\gamma(z))=(\iota\gamma)(z)=(\iota'(\gamma))(z)$ for all $z \in Z$. Hence $\iota'(\beta)\leq\iota'(\gamma)$.
\item  We show that for all $\beta,\gamma\in\Theta_{(Z,X)}, \beta\leq\gamma$, the diagram \[\begin{xy}\xymatrix{(J_{(Z,X)})_\beta\ar[r]^{\kappa'_\beta}\ar[d]_{(J_{(Z,X)})_{\beta\leq\gamma}}&(J_{(Z,Y)})_{\iota'(\beta)}\ar[d]^{(J_{(Z,Y)})_{\iota'(\beta)\leq\iota'(\gamma)}} \\(J_{(Z,X)})_{\gamma} \ar[r]_{\kappa'_{\gamma}} &(J_{(Z,Y)})_{\iota'(\gamma)}}\end{xy}\] commutes . Let $\beta,\gamma\in\Theta_{(Z,X)}, \beta\leq\gamma$, then \begin{align*}&\hspace{0.7cm}\kappa'_{\gamma}(J_{(Z,X)})_{\beta\leq\gamma}\\&=(\prod\limits_{z\in Z}\kappa_{\beta(z)})(\prod\limits_{z\in Z}G_{\beta(z)\leq\gamma(z)})\\ &\overset{(*)}{=}\prod\limits_{z\in Z}\kappa_{\beta(z)}G_{\beta(z)\leq\gamma(z)}\\& \overset{(**)}{=}\prod\limits_{z\in Z}G_{\iota(\beta(z))\leq\iota(\gamma(z))}\kappa_{\beta(z)}\\&=\prod\limits_{z\in Z}G_{\iota(\beta(z))\leq\iota(\gamma(z))}\prod\limits_{z\in Z}\kappa_{\beta(z)}\end{align*}\begin{align*}&=(J_{(Z,Y)})_{\iota\beta\leq\iota\gamma}\kappa'_\beta\\&=(J_{(Z,Y)})_{\iota'(\beta)\leq\iota'(\gamma)}\kappa'_\beta
\end{align*}.\\
(*) If $\sigma\in (J_{(Z,X)})_\beta$, then \begin{align*}&\hspace{0.4cm}\prod\limits_{z\in Z}\kappa_{\gamma(z)}(\prod\limits_{z\in Z}G_{\beta(z)\leq\gamma(z)}((\sigma_z)_{z\in Z})\\&=\prod\limits_{z\in Z}\kappa_{\gamma(z)}((G_{\beta(z)\leq\gamma(z)}(\sigma_z))_{z\in Z})\\&=(\kappa_{\gamma(z)}(G_{\beta(z)\leq\gamma(z)}(\sigma_z))_{z\in Z}\\&=\prod\limits_{z\in Z}\kappa_{\beta(z)}G_{\beta(z)\leq\gamma(z)}(\sigma_z)_{z\in Z}\end{align*}.\\
(**) The diagram \[\begin{xy}\xymatrix{G_{\beta(z)}\ar[r]^{\kappa_{\beta(z)}}\ar[d]_{G_{\beta(z)\leq\gamma(x)}}&H_{\iota(\beta(z))}\ar[d]^{H_{\iota(\beta(z))\leq\iota(\gamma(z))}} \\G_{\gamma(z)} \ar[r]_{\kappa_{\gamma(z)}} &H_{\iota(\gamma(z))}}\end{xy}\]
commutes for all $z\in Z$, since $(\iota,\kappa,\lambda)$ is a regular morphism.
\item We show that $Mor(1_Z,\lambda)(Mor(Z,X)_\beta)\subseteq Mor(Z,X)_{\iota'(\beta)}\ \forall\beta\in\Theta_{(Z,X)}$. Let $\beta\in\Theta_{(Z,X)}$, $f\in Mor(Z,X)_\beta$. By definition,  \begin{enumerate}[(1)]\item $f(z)\in X_{\beta(z)}\ \forall z\in Z$, and
\item if $\{z_0,...,z_p\}\subseteq Z$ is a local frame then $\{f(z_0),...,f(z_p)\}$ is a local b-frame in X for some $b\in\Phi$ such that $b\geq\beta(z_i)\ (0\leq i\leq p)$, i.e $\exists g_1,...,g_p\in G_b$ such that $g_i f(z_0)=f(z_i)\ (1\leq i\leq p)$.\end{enumerate} We must show that $\lambda f\in(Mor(Z,Y))_{\iota'(\beta)}$, i.e. \begin{enumerate}[(1')]
\item $(\lambda f)(z)\in X_{\iota(\beta(z))}\ \forall z\in Z$, and
\item if $\{z_0,...,z_p\}$ is a local frame in Z then $\{\lambda(f(z_0)),...,\lambda(f(z_p))\}$ is a local c-frame for some $c\in\Psi$ such that $c\geq\iota(\beta(z_i))\ (0\leq i\leq p)$, i.e. $\exists h_1,...,h_p\in H_c$ such that $ h_i\lambda(f(z_0))=\lambda(f(z_i))$.\end{enumerate} But (1) clearly implies (1'). Let now $\{z_0,...,z_p\}$ be a local frame in Z. Set $c=\iota(b)$ and $h_i=\kappa_b(g_i)\ (1\leq i\leq p)$.Then $h_i\lambda (f(z_0))=\kappa_b(g_i)\lambda(f(z_0))=\lambda(g_i f(z_0))=\lambda(f(z_i))\ (1\leq i\leq p)$ whence $\{z_0,...,z_p\}$ is a local c-frame in Y and obiously $c\geq\iota(\beta(z_i))\ (0\leq i\leq p)$, since $\iota$ is relation preserving.
\item We show that if $\beta\in\Theta_{(Z,X)}$ then $(\kappa'_\beta, Mor(1_Z,\lambda)|_{(Mor(Z,X))_\beta}):(J_{(Z,X)})_\beta\curvearrowright (Mor(Z,X))_\beta\rightarrow(J_{(Z,Y)})_{\iota(\beta)}\curvearrowright(Mor(Z,Y))_{\iota(\beta)}$ is a morphism of group actions, i.e. $Mor(1_Z,\lambda)(\sigma f)=\kappa'_\beta(\sigma)Mor(1_Z,\lambda)(f)\ \forall \sigma\in (J_{(Z,X)})_\beta, f\in(Mor(Z,X))_\beta$. Let $\beta\in\Theta_{(Z,X)}$, $\sigma\in (J_{(Z,X)})_\beta$,$f\in(Mor(Z,X))_\beta, z\in Z$. Then \begin{align*}&(Mor(1_Z,\lambda)(\sigma f))(z)=(\lambda(\sigma f))(z)=\lambda((\sigma f)(z))=\lambda(\sigma_z f(z))=\kappa_{\beta(z)}(\sigma_z)\lambda(f(z))\\&=(\kappa'_\beta(\sigma)(\lambda f))(z)=(\kappa'_\beta(\sigma)Mor(1_z,\lambda)(f))(z)\end{align*}
Hence $Mor(1_Z,\lambda)(\sigma f)=\kappa'(\sigma)Mor(1_Z,\lambda)(f)$. \hfill$\Box$\\
\end{enumerate}
}
\Corollary{A regular morphism is $\infty$-normal (more precisely if $(\iota,\kappa,\lambda)$ is a regular morphism then $\lambda$ is $\infty$-normal).}
\Proof{Let $(\iota,\kappa,\lambda):X\rightarrow Y$ be a regular morphism. By applying the lemma above repeatedly one can show that if $Z_1,...,Z_n$ is a sequence of global actions then $Mor(1_{Z_n},...,Mor(1_{Z_1},\lambda))...):Mor(Z_n,Mor(Z_{n-1},...,Mor(Z_1,X))...)\rightarrow Mor(Z_n,$ $ Mor(Z_{n-1}, ...,Mor(Z_1,Y))...)$ extends to a regular morphism and hence is a morphism of global actions (a formal proof would be one by induction). Thus $\lambda$ is $\infty$-normal.\hfill$\Box$} 
\Definition[product of global actions\index{global action!product of -s}]{
{Let $\{G_\alpha\curvearrowright X_\alpha|\alpha\in\Phi\}$ denote a global action on X and $\{H_\beta\curvearrowright Y_\beta|\beta\in\Psi\}$ a global action on Y. The global action on $X\times Y$ defined by $\Theta=\Phi\times \Psi$, $(\alpha,\beta)\leq(\alpha',\beta')\Leftrightarrow (\alpha\leq\alpha')\wedge(\beta\leq\beta')$, $J_{(\alpha,\beta)}=G_\alpha\times H_\beta$, $J_{(\alpha,\beta)\leq(\alpha',\beta')}=(G_{\alpha\leq\alpha'},H_{\beta\leq\beta'})$, $(X\times Y)_{(\alpha, \beta)}=X_\alpha\times X_\beta$, and $J_{(\alpha,\beta)}$ acts on $(X\times Y)_{(\alpha, \beta)}$ by $(g,h)(x,y)=(gx,hy)$ is called the \textit{product} of X and Y.
}\\

Our next goal is to show that the exponential law holds for many global action spaces. This will be needed for developing a good homotopy theory of global actions. 

The following notation will be used in the next definition. If A and B are sets, let $(A,B)=Mor_{((sets))}(A,B)$. If C is also a set then the exponential law for sets states that the function \begin{align}E_{set}:(A,(B,C))&\rightarrow (A\times B,C)\tag{2.18.1}\\ f&\mapsto E_{set}f\notag\end{align} where $E_{set}f(a,b)=f(a)(b) $ is an isomorphism of sets.

\Definition[exponential map\index{exponential map}]{Let $\{G_\alpha\curvearrowright X_\alpha|\alpha\in\Phi\}$ denote a global action on X ,  $\{H_\beta\curvearrowright Y_\beta|\beta\in\Psi\}$ a global action on Y and $\{J_\alpha\curvearrowright Z_\gamma|\gamma\in\Theta\}$ denote a global action on Z. Define a regular morphism $E=(\iota,\kappa,\lambda):Mor(X, Mor(Y,Z))\rightarrow Mor(X\times Y, Z)$ as follows. Define \[\iota:\begin{xy}\xymatrix{\Phi_{(X,(Y,Z))} \ar[r] \ar@{=}[d]&\Phi_{(X\times Y, Z)} \ar@{=}[d]\\(X, (Y,\Theta))&(X\times Y, \Theta)}\end{xy}\] to be the set theoretic exponential isomorphism (2.18.1). Define \[\kappa_\alpha:\begin{xy}\xymatrix{(J_{(X,(Y,Z))})_\alpha\ \ar[r] \ar@{=}[d]&(J_{(X\times Y, Z)})_{\iota(\alpha)} \ar@{=}[d]\\\prod \limits_{x\in X}(\prod \limits_{y\in Y}(J_{\alpha(x)(y)})&\prod\limits_{(x,y)\in X\times Y}J_{(\iota(\alpha))(x,y)}\ar@{=}[d]\\ & \prod\limits_{(x,y)\in X\times Y}J_{\alpha(x)(y)}}\end{xy} \] in the obvious way. Define $\lambda$ to be the composite mapping $Mor(X,Mor(Y,Z))\rightarrow(X,(Y,Z))\xrightarrow{E_{set}}(X\times Y,Z)$ where $E_{set}$ is the exponential morphism (2.18.1) (one can show that the image of the mapping defined above lies in $Mor(X \times Y,Z)$).}\\

We want to provide conditions when the map $Mor(X, Mor(Y,Z))\rightarrow Mor(X\times Y, Z)$  above is an isomorphism of global actions. This will be answered in Theorem 2.31 below.

\Definition[$\infty$-exponential\index{global action!$\infty$-exponential -}]{A global action Z is called \textit{$\infty$-exponential}, if for any global actions X and Y the regular morphism \[E:Mor(X, Mor(Y,Z))\rightarrow Mor(X\times Y, Z)\] above is an $\infty$-normal isomorphism. (This means  its inverse is $\infty$-normal). A global action Z is called \textit{regularly $\infty$-exponential}, if for any global actions X and Y, the regular morphism E above is a weak regular isomorphism. (This means its inverse is a regular morphism). It follows from Corollary 2.17 that Z is $\infty$-exponential also.
}\\

Suppose Z is $\infty$-exponential. Let $X_1,...,X_n$ be any n global actions where $n\geq2$. We construct now by induction on n the canonical $\infty$-normal isomorphism \[e_n:Mor(X_n,Mor(X_{n-1},...,Mor(X_1,Z))...)\rightarrow Mor(X_n\times...\times X_1,Z).\] Let $e_2:Mor(X_2$ $,Mor(X_1,Z)\rightarrow Mor(X_2\times X_1, Z)$ be the canonical (regular) morphism defined in 2.20. Since Z is $\infty$-exponential $e_2$ is by definition an $\infty$-normal isomorphism. Suppose $e_{n-1}:Mor(X_{n-1},...,Mor(X_1,Z))...)\rightarrow Mor(X_{n-1}\times...\times X_1,Z)$ has been constructed and is an  $\infty$-normal isomorphism. Define $e_n$ as the composition of the morphisms $(1_{X_n},e_{n-1}):Mor(X_n,Mor(X_{n-1},...,Mor(X_1,Z))...)\rightarrow Mor(X_n,Mor(X_{n-1}\times...\times X_1,Z))$, which is an $\infty$-normal isomorphism by the induction assumption for n-1 and Lemma 2.14, and the canonical morphism $Mor(X_n,$ $ Mor(X_{n-1}\times...\times X_1,Z))\rightarrow Mor(X_n\times...\times X_1,Z)$, which is an $\infty$-normal isomorphism because Z is $\infty$-exponential. As a composition of two $\infty$-normal isomorphisms,  $e_n$ is an $\infty$-normal isomorphism.\hfill$\Box$\\

The construction of $e_n$ just given yields the following lemma.
\Lemma{If Z is $\infty$-exponential then the canonical morphism \[e_n:Mor(X_n,Mor(X_{n-1},...,Mor(X_1,Z))...)\rightarrow Mor(X_n\times...\times X_1,Z)\] is an $\infty$-normal isomorphism.
}\\

\normalfont The pattern of the proof of Lemma 2.22 can be taken over to prove the following lemma.
\Lemma{If Z is regularly $\infty$-exponential then the canonical isomorphism $e_n$ defined in Lemma 2.22 extends to a weakly regular isomorphism.}
\Definition[{strong infimum condition}\index{global action!strong infimum -}]{Let $\{G_\alpha\curvearrowright X_\alpha|\alpha\in\Phi\}$ denote a global action on X. Let $\Delta\subseteq\Phi$ denote a finite subset and let $\Phi_{\geq\Delta}=\{\alpha\in\Phi|\alpha\geq\beta\ \forall\beta\in\Delta\}$. Let $U\subseteq X$ be a finite and nonempty subset such that $U\cap X_\beta\neq\emptyset$ for all $\beta\in\Delta$. The \textit{strong infimum condition} for X says that for any $\Delta$ and U as above, the set $\{\alpha\in\Phi_{\geq\Delta}|U\ \textnormal{is an $\alpha$-frame}\}$ is either empty or contains an initial element. If X satisfies the condition for $\Delta=\emptyset$ it is called an \textit{infimum action}\index{global action!infimum -}.
}\\

\Lemma{Let $\{G_\alpha\curvearrowright X_\alpha|\alpha\in\Phi\}$ denote a global action on X. If X satisfies the conditions (a) and (b) below then it is a strong infimum action.
\begin{enumerate}[(a)]
\item Let $\alpha,\beta\in\Phi$. Then $\alpha\leq\beta\Leftrightarrow\exists x\in X_\alpha\cap X_\beta$ such that $G_\alpha(x)\subseteq G_\beta(x)$. 
\item Let $\Psi\subseteq\Phi$. Then for any $x\in \bigcap\limits_{\alpha\in\Psi} X_\alpha$ there is a $\beta\in\Phi$ such that $\bigcap\limits_{\alpha\in\Psi}(G_\alpha(x))=G_\beta(x)$. 
\end{enumerate}}
\Proof{Let $\Delta\in\Phi$ denote a finite subset and $U\subseteq X$ denote finite and nonemtpy subset such that $U\cap X_\gamma\neq\emptyset\ \forall\gamma\in\Delta$. We must show that $\{\alpha\in\Phi_{\geq\Delta}|\textnormal{U is an }\alpha\textnormal{-frame}\}=:\Psi$ is either empty or contains an inital element. If $u\in U$ then clearly $u\in \bigcap\limits_{\alpha\in\Psi} X_\alpha$ and hence, by 2.25(b), there is a $\beta\in\Phi$ such that $\bigcap\limits_{\alpha\in\Psi}(G_\alpha(u))=G_\beta(u)$. Since $U\subseteq\bigcap\limits_{\alpha\in\Psi}(G_\alpha(u))=G_\beta(u)$ it is obvious that U is a $\beta$-frame. Since $G_\beta(u)\subseteq G_\alpha(u)\ \forall \alpha\in\Psi$, it follows by 2.25(a) that $\beta\leq\alpha\ \forall\alpha\in\Psi$. It remains to show that $\gamma\leq\beta\ \forall \gamma\in\Delta$. Let  $\gamma\in\Delta$. Since $U\cap X_{\gamma'}\neq\emptyset\ \forall\gamma'\in\Delta$ there is an $x\in U\cap X_\gamma$. Since $\gamma\leq\alpha\ \forall\alpha\in\Psi$ and $U\subseteq X_\alpha\ \forall\alpha\in\Psi$ it follows from 2.1(a) that $G_\gamma(x)\subseteq\bigcap\limits_{\alpha\in\Psi}(G_\alpha(x))=G_\beta(x)$. Hence, by 2.25(b), $\gamma\leq\beta$.\hfill$\Box$
}
\Corollary{A standard single domain global action is a strong infimum action
}
\Proof{We show that a standard single domain global action satisfies 2.25(a) and (b). \begin{enumerate}[(a)]
\item Let $\alpha,\beta\in\Phi$.\\
$\Rightarrow$: Suppose $\alpha\leq\beta$. We have to show that $\exists x\in X_\alpha\cap X_\beta=G$ such that $G_\alpha(x)\subseteq G_\beta(x)$. $\alpha\leq\beta$ implies $G_\alpha\subseteq G_\beta$. Hence $G_\alpha(1)=G_\alpha\subseteq G_\beta=G_\beta(1)$ where 1 denotes the identity element in G.\\
$\Leftarrow$: Suppose $\exists x\in X_\alpha\cap X_\beta=G$ such that $G_\alpha(x)\subseteq G_\beta(x)$. Let $g\in G_\alpha$. Then $gx\in G_\alpha(x)$. Since $G_\alpha(x)\subseteq G_\beta(x)$, $\exists g'\in G_\beta$ such that $gx=g'x$. Since $G_\alpha$ and $G_\beta$ act on $G$ by left or right multiplication, it follows that $g=g'$. Thus $G_\alpha\subseteq G_\beta$ and hence $\alpha\leq\beta$.
\item Let $\Psi\subseteq\Phi$ be a nonemtpy subset and $x\in G$. We have to show that there is a $\beta\in\Phi$ such that $\bigcap\limits_{\alpha\in\Psi}G_\alpha(x)=G_\beta(x)$. Since $\{G_\gamma|\gamma\in\Phi\}$ is closed under taking intersections, $\exists\beta\in\Phi$ such that $\bigcap\limits_{\alpha\in\Psi}G_\alpha=G_\beta$. Thus $(\bigcap\limits_{\alpha\in\Psi}G_\alpha)(x)=G_\beta(x)$. We show now that $(\bigcap\limits_{\alpha\in\Psi}G_\alpha)(x)=\bigcap\limits_{\alpha\in\Psi}(G_\alpha(x))$.\\
$\subseteq$: Let $y\in(\bigcap\limits_{\alpha\in\Psi}G_\alpha)(x)$. Then $\exists g\in(\bigcap\limits_{\alpha\in\Psi}G_\alpha)$ such that $y=gx$. $g\in(\bigcap\limits_{\alpha\in\Psi}G_\alpha)$ implies $g\in G_\alpha\ \forall\alpha\in\Psi$. It follows that $gx\in G_\alpha x\ \forall\alpha\in\Psi$ and hence $y=gx\in\bigcap\limits_{\alpha\in\Psi}(G_\alpha(x))$.\\
$\supseteq$: Let $z\in\bigcap\limits_{\alpha\in\Psi}(G_\alpha(x))$. Then $z\in G_\alpha(x)\ \forall\alpha\in\Psi$. Hence $\forall\alpha\in\Psi\ \exists g_\alpha\in G_\alpha$ such that $y=g_\alpha x$. If $\alpha,\beta\in\Psi$ then $g_\alpha x=y=g_\beta x$. But this implies $g_\alpha=g_\beta$. Hence if $\beta\in\Psi$ then $g_\beta\in G_\alpha\ \forall\alpha\in\Psi$. Thus $g_\beta\in\bigcap\limits_{\alpha\in\Psi}G_\beta$ and hence $z=g_\beta x\in(\bigcap\limits_{\alpha\in\Psi}G_\alpha)(x)$. \hfill$\Box$
\end{enumerate}
}

The four lemmas below will be used in the proof of Theorem 2.31. Their proofs will be given after that of Theorem 2.31.
\Lemma[Local-Global Lemma\index{Local-Global Lemma}]{Let X,Y denote global actions and $f_0\in Mor(X,Y)_\beta$. Then $\{f_0,...,f_p\}$ is a $\beta$-frame in $Mor(X,Y)$ if and only if $\{f_0(x),$ $f_1(x),...,f_p(x)\}$ is a $\beta(x)$-frame in Y \hspace{0.07cm}$\forall x \in X$.
}
\Lemma{Let X and Y be global actions and let $\{f_0,...,f_p\}$ be a $\beta$-frame in $Mor(X,Y)$. If $\{x_0,...x_q\}$ is a local frame in X then $\{f_i(x_j)|0\leq i\leq p, 0\leq j\leq q\}$ is a  b-frame in Y for some b such that $b\geq\beta(x_0),...,\beta(x_q)$.
}
\Lemma{An infimum action is conormal}.
\Lemma{If Z is an infimum action then for any global action X, $Mor(X,Z)$ is an infimum action. If Z is an strong infimum action and the relation on $\Theta$ is transitive then for any global action X, $Mor(X,Z)$ is a strong infimum action and the relation on $\Theta_{(X,Z)}$ is transitive.}
\Theorem[exponential law\index{exponential law}]{An infimum action is $\infty$-conormal and $\infty$-exponential. A strong infimum action is $\infty$-conormal and regularly $\infty$-exponential.}
\Proof{Let Z be an infimum action.\\
\underline{Z is $\infty$-conormal:} We have to show that if Y is a global action and $g:Y\rightarrow Z$ is a morphism then g is $\infty$-normal, i.e. if $X_n,...,X_1$ are global actions then the map \begin{eqnarray*} & & Mor(1_{X_n},Mor(1_{X_{n-1}},...,Mor(1_{X_1},g))...):\\& & Mor(X_n,Mor(X_{n-1},...,Mor(X_1,Y))...)\rightarrow Mor(X_n,Mor(X_{n-1},...,Mor(X_1,Z))...) \end{eqnarray*} is a morphism of global actions. We prove this by induction on n.\\
n=1: Let X be a global action. By Lemma 2.29, Z is conormal and hence the map $Mor(1_X,g):Mor(X,Y)\rightarrow Mor(X,Z)$ is a morphism of global actions.\\
$n-1\rightarrow n$: Let $n\geq2$. Suppose by induction that for any global actions $X'_{n-1},...,X'_1$, the map \begin{eqnarray*} & & Mor(1_{X'_{n-1}},Mor(1_{X'_{n-1}},...,Mor(1_{X'_1},g))...):\\& & Mor(X'_{n-1},...,Mor(X'_1,Y)...)\rightarrow Mor(X'_{n-1},...,Mor(X'_1,Z)...) \end{eqnarray*} is a morphism of global actions. Let $X_n,...X_1$ be global actions. By Lemmas 2.29 and 2.30, $Mor(X_{n-1},...,Mor(X_1,Z)...)$ is conormal and since by the induction assumption $Mor(1_{X_{n-1}},Mor(1_{X_{n-1}},...,Mor(1_{X_1},g))...)$ is a morphism of global actions, it follows that \begin{eqnarray*} & & Mor(1_{X_n},Mor(1_{X_{n-1}},...,Mor(1_{X_1},g))...):\\& & Mor(X_n,Mor(X_{n-1},...,Mor(X_1,Y))...)\rightarrow Mor(X_n,Mor(X_{n-1},...,Mor(X_1,Z))...) \end{eqnarray*} is a morphism of global actions.\\
\underline{Z is $\infty$-exponential:} We have to show that if X and Y are global actions then the morphism $E:Mor(X,$ $Mor(Y,Z))\rightarrow Mor(X\times Y,Z)$ is an $\infty$-normal isomorphism. Since E is regular, it is $\infty$-normal by Corollary 2.18. Therefore it sufficies to show that E has an $\infty$-normal inverse. But $Mor(X,Mor(Y,Z))$ is an infimum action by 2 applications of Lemma 2.30 and therefore $\infty$-conormal. Thus any inverse for E is an $\infty$-normal inverse.\\

Let X and Y be global actions. Define \begin{align*}E':Mor(X\times Y,Z)&\rightarrow Mor(X, Mor(Y,Z))\\ h&\mapsto E'(h) \end{align*}
where \begin{align*}E'(h):X&\rightarrow Mor(Y,Z)\\ x&\mapsto E'(h)(x) \end{align*} and \begin{align*}E'(h)(x):Y&\rightarrow Z\\ y&\mapsto h(x,y)\tag*{.}\end{align*} Then $E'$ is a set theoretical inverse of E. We want to show that $E'$ is a morphism of global actions. There are a number of things to verify.\\
\\
$\bullet$ $E'(h)(x):Y\rightarrow Z$ is a morphism of global actions $\forall h \in Mor(X\times Y, Z)$ and $\forall x\in X$.\\
Let $h \in Mor(X\times Y, Z)$ and $x\in X$. Let $\{y_0,...,y_p\}$ be a local frame in Y. We have to show that $\{E'(h)(x)(y_0),...,$ $E'(h)(x)(y_p)\}$ is a local frame in Z. Since $X=\bigcup\limits_{\alpha\in\Phi}X_\alpha$ there is an $\alpha\in\Phi$ such that $x\in X_\alpha$. Since $\{y_0,...,y_p\}$ is a local frame there is an $\beta\in\Psi$ such that $\{y_0,...,y_p\}\subseteq Y_\beta$ and is a $\beta$-frame. Clearly $\{(x,y_0),...,(x,y_p)\}\subseteq X_\alpha\times Y_\beta=(X\times Y)_{\alpha,\beta}$.  Since $\{y_0,...,y_p\}$ is a  $\beta$-frame there are $h_1,...,h_p\in H_\beta$ such that $h_iy_0=y_i$ $(1\leq i\leq p)$. 
If 1 denotes the identity element in the local group $G_\alpha$ then  $(1,h_i)(x,y_0)=(1x,h_iy_0)=(x,y_i)$ $(1\leq i\leq p)$. Thus $\{(x,y_0),...,(x,y_p)\}$ is a local frame and since $h:X\times Y\rightarrow Z$ is a morphism of global actions $\{h(x,y_0),...,h(x,y_p)\}=\{E'(h)(x)(y_0),...,E'(h)(x)(y_p)\}$ is a local frame in Z. Thus $E'(h)(x)$ is a morphism of global actions.\\
\\
$\bullet$ $E'(h):X\rightarrow Mor(Y,Z)$ is a morphism of global actions $\forall h \in Mor(X\times Y,Z)$.\\Let $h \in Mor(X\times Y,Z)$ and $\{x_0,...,x_p\}$ be a local frame in X. We must show that $\{E'(h)(x_0),...,E'(h)(x_p)\}$ is a local frame in $Mor(Y,Z)$. As above $\{(x_0,y),...,(x_p,y)\}$ is a local frame in  $X\times Y$ for any $y\in Y$. It follows that $\{h(x_0,y),...,h(x_p,y)\}=\{E'(h)(x_0)(y),$ $...,E'(h)(x_p)(y)\}$ is a local frame in Z \hspace{0cm}for any $y\in Y$ because $h:X\times Y\rightarrow Z$ is a morphism of global actions. Since Z is an infimum action, the set $\{c\in\Theta|\{E'(h)(x_0)(y),$ $...,E'(h)(x_p)(y)\}\ \textnormal{is a c-frame}\}$ has an initial element $c_y$ for any $y\in Y$. Define \begin{align*}\gamma:Y&\rightarrow\Theta\\ y&\mapsto c_y\ .\end{align*} We shall show that $E'(h)(x_0)\in (Mor(Y,Z))_\gamma$, i.e. \begin {enumerate} [(a)]
\item $E'(h)(x_0)(y)\in Z_{\gamma(y)}\ \forall y\in Y$, and
\item if $\{y_0,...,y_q\}$ is a local frame in Y then $\{E'(h)(x_0)(y_0),...,E'(h)(x_0)(y_p)\}$ is a c-frame in Z for some  $c\geq \gamma(y_i)\ (0\leq i\leq q)$. 
\end{enumerate}
Clearly if $y\in Y$ then $\{E'(h)(x_0)(y),$ $...,E'(h)(x_p)(y)\}$ is a $c_y$-frame. Thus $E'(h)$ $(x_0)(y)\in Z_{c_y}=Z_{\gamma(y)}\ \forall y\in Y$. (b) holds because of the following. Let $\{y_0,...,y_q\}$ be a local frame in Y. Then $\{(x_i,y_j)|0\leq i\leq p, 0\leq j\leq q\}$ is a local frame in $X\times Y$ and hence $\{h(x_i,y_j)|0\leq i\leq p, 0\leq j\leq q\}$=$\{E'(h)(x_i)(y_j)|0\leq i\leq p, 0\leq j\leq q\}$ is a c-frame for some $c\in\Theta$. On the one hand this implies that $\{E'(h)(x_0)(y_0),...,E'(h)(x_0)(y_p)\}$ is a c-frame and on the other hand that $\{E'(h)(x_0)(y_j),...,E'(h)(x_p)(y_j)\}$ is a c-frame $\forall j\in\{0,...,q\}$ and hence $c\geq c_{y_j}=\gamma(y_j)\ \forall j\in\{0,...,q\}$. Thus $E'(h)(x_0)\in (Mor(Y,Z))_\gamma$. By the Local-Global Lemma 2.27 it follows that $\{E'(h)(x_0),...,E'(h)(x_p)\}$ is a local frame.\\
\\
$\bullet$ $E':Mor(X\times Y,Z))\rightarrow Mor(X, Mor(Y,Z))$ is a morphism of global actions.\\
Let $\{h_0,...,h_p\}$ be a local frame in $Mor(X\times Y,Z)$. We have to show that $\{E'(h_0),...,$ $E'(h_p)\}$ is a local frame in $Mor(X, Mor(Y,Z))$. By Lemma 2.27, $\{h_0(x,y),...,$ $h_p(x,y)\}=\{E'(h_0)(x)(y),...,E'(h_p)(x)(y)\}$ is a local frame in Z  $\forall(x,y)\in X\times Y$. Since Z is an infimum action, the set $\{c\in\Theta|\{h_0(x,y),...,h_p(x,y)\} \textnormal{ is a c-frame}\}$ has an initial element $c_{(x,y)}$. Define \begin{align*}\gamma:X&\rightarrow(Y,\Theta)\\ x&\mapsto c_{(x,-)}\ .\end{align*} We will now show that $E'(h_0)(x)\in Mor(Y,Z)_{\gamma(x)}$ for any $x\in X$, i.e. \begin{enumerate}[(a)]
\item $E'(h_0)(x)(y)=h_0(x,y)\in Z_{\gamma(x)(y)=c(x,y)}\ \forall y\in Y$, and
\item if $\{y_0,...,y_r\}$ is a local frame in Y then $\{E'(h_0)(x)(y_0),...,E'(h_0)(x)(y_r)\}=\{h_0(x,y_0),...,h_0(x,y_r)\}$ is a d-frame in Z for some $d\geq\gamma(x)(y_j)=c(x,y_j)\ (0\leq j\leq r)$.
\end{enumerate} 
Clearly $\{h_0(x,y),...,h_p(x,y)\}$ is a $c(x,y)$-frame $\forall (x,y)\in X\times Y$. Thus $h_0(x,y)\in Z_{c(x,y)}\ \forall (x,y)\in X\times Y$. (b) holds because of the following. Let $\{y_0,...,y_r\}$ be a local frame in Y. Then $\{(x,y_0),...,(x,y_r)\}$ is a local frame in $X\times Y$. By Lemma 2.28, $\{h_i(x,y_j)|0\leq i\leq p, 0\leq j\leq r\}$ is a d-frame in Z. But on the one hand this implies that $\{h_0(x,y_0),...,h_0(x,y_r)\}$ is a d-frame and on the other hand that $\{h_0(x,y_j),...,h_p(x,y_j)\}$ is a d-frame for any $j\in\{0,...,r\}$ and hence $e\geq c(x,y_j)\ (0\leq j \leq r)$. Thus (b) holds and hence $E'(h_0)(x)\in Mor(Y,Z)_{\gamma(x)}$ for any $x\in X$. Since $\{E'(h_0)(x)(y),...,$ $E'(h_p)(x)(y)\}$ is a $\gamma(x)(y)$-frame $\forall x\in X$ and $\forall y\in Y$ it follows from Lemma 2.27 that $\{E'(h_0)(x),...,E'(h_p)(x)\}$ is a $\gamma(x)$-frame in $Mor(Y,Z)\ \forall x\in X$.\\

We show now that that $E'(h_0)\in (Mor(X,Mor(Y,Z)))_\gamma$, i.e. \begin{enumerate} [(1)]
\item $E'(h_0)(x)\in Mor(Y,Z)_{\gamma(x)}\ \forall x\in X$,  and
\item if $\{x_0,...,x_q\}$ is a local frame in X then $\{E'(h_0)(x_0),...,E'(h_0)(x_q)\}$ is a c-frame in $Mor(Y,Z)$ for some $c\geq\gamma(x_i)\ (0\leq i\leq q)$.
\end{enumerate} 
We have shown (1) above. Let $\{x_0,...,x_q\}$ be a local frame in X. Then $\{(x_0,y),...,$ $(x_q,y)\}$ is a local frame in $X\times Y$ for any $y\in Y$. It follows from Lemma 2.28 that $\{h_i(x_j,y)|0\leq i\leq p, 0\leq j\leq q\}$ is a local frame in Z for any $y \in Y$. Since Z is an infimum action, for any $y \in Y$ the set $\{d\in\Theta|\{h_i(x_j,y)|0\leq i\leq p, 0\leq j\leq q\}\textnormal{ is a local frame}\}$ has an initial element $d_y$.  Define \begin{align*}\delta:Y&\rightarrow\Theta\\ y&\mapsto d_y\ .\end{align*} Since $\{h_i(x_j,y)|0\leq i\leq p, 0\leq j\leq q\}$ is a $\delta(y)$-frame for any $y\in Y$, it follows that $\{h_0(x_j,y)|0\leq j\leq q\}$=$\{E'(h_0)(x_j)(y)|0\leq j\leq q\}$ is a $\delta(y)$-frame for any $y\in Y$. We show now that $E'(h)(x_0)\in Mor(Y,Z)_\delta$, i.e. \begin{enumerate}[({1.}1)] 
\item $E'(h_0)(x_0)(y)=h_0(x_0,y)\in Z_{\delta(y)=d_y}\ \forall y\in Y$, and
\item if $\{y_0,...,y_r\}$ is an local frame in Y then $\{E'(h_0)(x_0)(y_0),...,E'(h_0)(x_0)(y_r)\}=\{h_0(x_0,y_0),...,h_0(x_0,y_r)\}$ is an e-frame in Z for some $e\geq\delta(y_j)=d_{y_j}\ (0\leq j\leq r)$.
\end{enumerate} Clearly if $y\in Y$ then $\{h_0(x_j,y)|0\leq j\leq q\}$ is a $\delta(y)$-frame. Thus $h_0(x_0,y)\in Z_{\gamma(x)}\ \forall y\in Y$. We show that (1.2) holds. Let $\{y_0,...,y_r\}$ be a local frame in Y. Since $\{x_j,y_k|0\leq j\leq q, 0\leq k\leq r\}$ is a local frame in $X\times Y$, it follows from Lemma 2.28 that $\{h_i(x_j,y_k)|0\leq i\leq p, 0\leq j\leq q, 0\leq k\leq r\}$ is an e-frame in Z. On the one hand this implies that $\{h_0(x_0,y_0),...,h_0(x_0,y_r)\}$ is an e-frame and on the other hand that $\{h_i(x_j,y_k)|0\leq i\leq p, 0\leq j\leq q\}$ is an e-frame for any $k\in\{0,...,r\}$ and hence $e\geq d_{y_j}=\delta(y_j)\ (0\leq k\leq r)$. Thus $E'(h)(x_0)\in Mor(Y,Z)_\delta$. By Lemma 2.27, $\{E'(h_0)(x_j)|0\leq j\leq q\}$ is a $\delta$-frame. Obviously if $y\in Y$ then $d_y\geq c_{(x_j,y)}\ (0 \leq j\leq q)$ and hence $\delta\geq\gamma(x_j)\ (0\leq j\leq q)$. Thus (2) holds and hence $E'(h_0)\in (Mor(X,Mor(Y,Z)))_\gamma$. Since we have shown that $\{E'(h_0)(x),...,E'(h_p)(x)\}$ is a $\gamma(x)$-frame in $Mor(Y,Z)\ \forall x\in X$ it follows from Lemma 2.27 that $\{E'(h_0),...,$ $E'(h_p)\}$ is a $\gamma$-frame.\\ 

Let Z now be a strong infimum action. We must show that Z is regularly $\infty$-exponential, i.e. if X and Y are global actions then $E=(\iota,\kappa,\lambda):Mor(X,(Mor(Y$ $,Z))\rightarrow Mor(X\times Y, Z)$ is a weak regular isomorphism, i.e. there is a regular morphism $(\iota',\kappa',\lambda'):Mor(X\times Y, Z)\rightarrow Mor(X,(Mor(Y,Z))$ such that $\lambda'$ is inverse to $\lambda$. Let X and Y be global actions. Define $\iota'$ as the set theoretic inverse of $\iota$. Define  $\lambda'=E'$ where $E'$ is the morphism of global actions we constructed above. Define \[\kappa'_\alpha:\begin{xy}\xymatrix{(J_{(X\times Y,Z)})_\alpha\ \ar[r] \ar@{=}[d]&(J_{(X,Mor(Y,Z))})_{\iota'(\alpha)} \ar@{=}[d]\\ \prod \limits_{(x,y)\in X\times Y}J_{\alpha(x,y)}&\prod\limits_{x\in X}(\prod\limits_{y\in Y}(J_{(\iota'(\alpha))(x,y)})\ar@{=}[d]\\ & \prod\limits_{x\in X}(\prod\limits_{y\in Y}(J_{\alpha(x,y)})}\end{xy}\] in the obvious way. Since  $\lambda'=E'$ is inverse to $\lambda$ it only remains to show that $(\iota',\kappa',\lambda')$ is a regular morphism, i.e. \begin {enumerate}[(a)]
\item $\iota'(\alpha)\leq\iota'(\beta)$ whenever $\alpha\leq\beta$, 
\item if $\alpha\leq\beta$ then the diagram \[\begin{xy}\xymatrix{(J_{(X\times Y, Z)})_\alpha \ar[r]^{\kappa'_\alpha} \ar[d]_{(J_{(X\times Y, Z)})_{\alpha\leq\beta}}&(J_{(X,(Y,Z))})_{\iota'(\alpha)} \ar[d]^{(J_{(X,(Y,Z))})_{\iota'(\alpha)\leq\iota'(\beta)}} \\(J_{(X\times Y, Z)})_\beta \ar[r]_{\kappa'_\beta} &(J_{(X,(Y,Z))})_{\iota'(\beta)}}\end{xy}\] commutes,  
\item $\lambda'((Mor(X\times Y,Z))_\alpha)\subseteq (Mor(X,Mor(Y,Z)))_{\iota'(\alpha)}\ \forall\alpha\in\Theta_{(Mor(X\times Y,Z))}$, 
\item for all $\alpha\in\Theta_{((X\times Y,Z))}$ the pair $(\kappa'_\alpha, \lambda'|(Mor(X\times Y,Z)_\alpha):(J_{(X\times Y,Z))})_\alpha\curvearrowright (Mor(X\times Y,Z))_\alpha\rightarrow (J_{(Mor(X,Mor(Y,Z))})_{\iota'(\alpha)} \curvearrowright (Mor(X,Mor(Y,Z))_{\iota'(\alpha)}$ is a morphism of group actions.
\end{enumerate} We will now verify (a)-(d).\begin{enumerate}[(a)]
\item Let $\alpha,\beta\in\Theta_{(X\times Y,Z)}=(X\times Y,\Theta)$ such that $\alpha\leq\beta$. We have to show that $\iota'(\alpha)\leq\iota'(\beta)$. But \begin{align*}&\hspace{0.64cm}\iota'(\alpha)\leq\iota'(\beta)\\ &\Leftrightarrow\iota'(\alpha)(x)\leq\iota'(\beta)(x)\ \forall x\in X\\ &\Leftrightarrow \iota'(\alpha)(x)(y)\leq\iota'(\beta)(x)(y)\ \forall x\in X, y\in Y\\ &\Leftrightarrow \alpha(x,y)\leq\beta(x,y). \end{align*}
\item Let $\alpha,\beta\in\Theta_{(X\times Y,Z)}=(X\times Y,\Theta)$ such that $\alpha\leq\beta$ and $\sigma=(\sigma_{(x,y)})_{(x,y)}\in(J_{(X\times Y,Z)})_\alpha$. Then \begin{align*}&\hspace{0.55cm}\kappa'_\beta((J_{(X\times Y,Z)})_{\alpha\leq\beta}(\sigma))\\&=\kappa'_\beta(\prod\limits_{(x,y)\in X\times Y}J_{\alpha(x,y)\leq\beta(x,y)}(\sigma))\\&=\kappa'_\beta((J_{\alpha(x,y)\leq\beta(x,y)}(\sigma_{(x,y)}))_{(x,y)})\\ &=((J_{\alpha(x,y)\leq\beta(x,y)}(\sigma_{(x,y)}))_y)_x\\ &=((J_{\iota'(\alpha)(x)(y)\leq\iota'(\beta)(x)(y)}(\sigma_{(x,y)}))_y)_x\\ &=(\prod\limits_{x\in X}(\prod\limits_{y\in Y}J_{\iota'(\alpha)(x)(y)\leq\iota'(\beta)(x)(y)}))(((\sigma_{(x,y)})_y)_x)\\ &=(\prod\limits_{x\in X}(J_{(Y,Z)})_{\iota'(\alpha)(x)\leq\iota'(\beta)(x)})(((\sigma_{(x,y)})_y)_x)\\ &=(J_{(X,(Y,Z))})_{\iota'(\alpha)\leq\iota'(\beta)}(((\sigma_{(x,y)})_y)_x)\\ &=(J_{(X,(Y,Z))})_{\iota'(\alpha)\leq\iota'(\beta)}(\kappa'_\alpha(\sigma))\ \textnormal{,} \end{align*}
i.e. $\kappa'_\beta\circ (J_{(X\times Y,Z)})_{\alpha\leq\beta}=(J_{(X,(Y,Z))})_{\iota'(\alpha)\leq\iota'(\beta)}\circ\kappa'_\alpha$.
\item Let $\alpha\in\Theta_{(X\times Y,Z)}$ and $f\in(Mor(X\times Y,Z))_\alpha$, i.e. \begin{enumerate}[(1)]
\item $f(x,y)\in Z_{\alpha(x,y)}\ \forall(x,y)\in X\times Y$, and
\item if $\{(x_0,y_0),...,(x_p,y_p)\}$ is a local frame in $X\times Y$ then $\{f(x_0,y_0),...,f(x_p,y_p)\}$ is a c-frame for some $c\in\Theta$ such that $c\geq\alpha(x_i,y_i)\ (0\leq i\leq p)$.\end{enumerate} We have to show that  $\lambda'(f)\in(Mor(X,Mor(Y,Z)))_{\iota'(\alpha)}$, i.e. \begin{enumerate}[(1$'$)]
\item $\lambda'(f)(x)\in ((Mor(Y,Z))_{\iota'(\alpha)(x)}\ \forall x\in X$, and
\item if $\{x_0,...,x_p\}$ is a local frame in X then $\{\lambda'(f)(x_0),...,\lambda'(f)(x_p)\}$ is a d-frame for some $d\in\Theta_{(Y,Z)}$ such that $d\geq\iota'(\alpha)(x_i)\ (0\leq i\leq p)$.\end{enumerate} First we show (1$'$), i.e. if $x\in X$ then \begin{enumerate}[(1$''$)]
\item $\lambda'(f)(x)(y)=f(x,y)\in Z_{\iota'(\alpha)(x)(y)=\alpha(x,y)}\ \forall y\in Y$, and
\item if $\{y_0,...,y_q\}$ is a local frame in Y then $\{\lambda'(f)(x)(y_0),...,\lambda'(f)(x)(y_q)\}=\{f(x,y_0),$ $...,f(x,y_q)\}$ is an e-frame for some $e\in\Theta$ such that $e\geq\iota'(\alpha)(x)(y_j)=\alpha(x,y_j)\ (0\leq j\leq q)$.\end{enumerate} (1$''$) follows from (1). Let $\{y_0,..,y_q\}$ be a local frame in Y. Then $\{(x,y_0),...,(x,y_q)\}$ is a local frame in $X\times Y$ and hence it follows from (2) that $\{f(x,y_0),...,f(x,y_q)\}$ is an e-frame for some $e\in\Theta$ such that $e\geq\alpha(x,y_q)\ (0\leq j\leq q)$. Thus (2$''$) holds and hence (1$'$) holds.\\

We show now (2$'$). Let $\{x_0,...,x_p\}$ be a local frame in X. Then for any $y\in Y$ $\{(x_0,y),...,(x_p,y)\}$ is a local frame in $X\times Y$ and since $f:X\times Y\rightarrow Z$ is a morphism of global actions $\{f(x_0,y),...,f(x_p,y)\}$ is a local frame in Z for any $y\in Y$. Since Z is a strong infimum action the set $\{c\in Z|\{f(x_0,y),...,f(x_p,y)\}\textnormal{ is a}$ $\textnormal{ c-frame},  c\geq\iota'(\alpha)(x_i)(y)\ (0\leq i\leq p)\}$ has an initial element $c_y$. Define \begin{align*}\gamma:Y&\rightarrow\Theta\\ y&\mapsto c_y.\end{align*} Since $\gamma(y)=c_y\geq\iota'(\alpha)(x_i)(y)\ \forall y\in Y, i\in\{0,...,p\}$ it follows that $\gamma\geq\iota'(\alpha)(x_i)\ (0\leq i\leq p)$. We show now that $\lambda'(f)(x_0)\in Mor(Y,Z)_\gamma$, i.e. \begin{enumerate}[(1$'''$)]
\item $\lambda'(f)(x_0)(y)=f(x_0,y)\in Z_{\gamma(y)=c_y}\ \forall y\in Y$, and
\item if $\{y_0,...,y_q\}$ is a local frame in Y then $\{\lambda'(f)(x_0)(y_0),...,\lambda'(f)(x_0)(y_q)\}=\{f(x_0,y_0),...,f(x_0,y_q)\}$ is an e-frame for some $e\in\Theta$ such that $e\geq\gamma(y_j)=c_{y_j}\ (0\leq j\leq q)$.\end{enumerate} Clearly $\{f(x_0,y),...,f(x_p,y)\}$ is a $c_y$-frame for any $y\in Y$. Thus $f(x_0,y)\in Z_{c_y}\ \forall y\in Y$. (2$'''$) holds because of the following. Let $\{y_0,...,y_q\}$ be a local frame in Y. Then $\{(x_i,y_j)|0\leq i \leq p, 0\leq j\leq q\}$ is a local frame in $X\times Y$ and hence it follows from (2) that $\{f(x_i,y_j)|0\leq i \leq p, 0\leq j\leq q\}$ is a c-frame for some $c\in\Theta$ such that $c\geq\alpha(x_i,y_j)\ (0\leq i\leq p, 0\leq j\leq q)$. This implies that $\{f(x_0,y_0),...,f(x_0,y_q)\}$ is a c-frame and that $\{f(x_0,y_j),...,f(x_p,y_j)\}$ is a c-frame for any $j\in\{0,...,q\}$, whence $c\geq c_{y_j}\ (0\leq j\leq q)$. Thus (2$'''$) holds and hence $\lambda'(f)(x_0)\in Mor(Y,Z)_\gamma$. Since $\{\lambda'(f)(x_0)(y),...,\lambda'(f)(x_p)(y)\}$=$\{f(x_0,y),$ $...,f(x_p,y)\}$ is a $c_y$-frame for any $y\in Y$, it follows from Lemma 2.27 that $\{\lambda'(f)(x_0),...,\lambda'(f)(x_p)\}$ is a $\gamma$-frame.
\item Let $\alpha\in\Theta_{(X\times Y,Z)}$, $\sigma\in(J_{X\times Y,Z)})_\alpha$ and $f\in Mor(X\times Y,Z)_\alpha$. We must show that $\lambda'(\sigma f)=\kappa'_\alpha(\sigma)\lambda'(f)$. Let $x\in X, y\in Y$. Then \begin{align*}&\hspace{0.5cm}\lambda'(\sigma f)(x)(y)\\&=\sigma_{(x,y)}f(x,y)\\&=\sigma_{(x,y)}\lambda'(f)(x)(y)\\&=((\sigma_{(x,y)})_y)_x\lambda'(f)\\&=\kappa'_\alpha(\sigma)\lambda'(f)\end{align*} and hence $\lambda'(\sigma f)=\kappa'_\alpha(\sigma)\lambda'(f)$. \hfill$\Box$\end{enumerate}
}
\Proof{\textbf{of Lemma 2.27} \hspace{0,1cm}Let $X,Y$ denote global actions and $f_0\in Mor(X,Y)_\beta$.\\ $\Rightarrow$: Let $\{f_0,...f_p\}$ be a $\beta$-frame in $Mor(X,Y)$, i.e. $\{f_0,...,f_p\}\subseteq Mor(X,Y)_\beta$ and $\exists\sigma_1,...,\sigma_p$ $\in J_\beta$ such that $\sigma_i f_0=f_i\ (1\leq i\leq p)$. We must show that $\{f_0(x),...,f_p(x)\}$ is a $\beta(x)$-frame $\forall x\in X$, i.e. if $x\in X$ then $\{f_0(x),...,f_p(x)\}\subseteq Y_{\beta(x)}$ and $\exists h_1,...h_p\in H_\beta(x)$ such that $h_if_0(x)=f_i(x)\ (1\leq i\leq p)$. But if $i\in\{0,...,p\}$ then $\sigma_i f_0=f_i\Leftrightarrow (\sigma_i f_0)(x)=f_i(x)\ \forall x\in X\Leftrightarrow (\sigma_i)_x f_0(x)=f_i(x)\ \forall x\in X$ and from $\{f_0,...,f_p\}\subseteq Mor(X,Y)_\beta$ it follows that $\{f_0(x),...,f_p(x)\}\subseteq Y_{\beta(x)}\ \forall x\in X$.\\
$\Leftarrow$: Let $\{f_0(x),...f_p(x)\}$ be a $\beta(x)$-frame in Y $\forall x\in X$, i.e. if $x\in X$ then $\{f_0(x),...,f_p(x)\}\subseteq Y_{\beta(x)}$ and $\exists h_1(x),...,h_p(x)\in H_{\beta(x)}$ such that $h_i(x)f_0(x)=f_i(x)\ (1\leq i\leq p)$. We have to show that $\{f_0,...,f_p\}$ is a $\beta$-frame, i.e. $\{f_0,...,f_p\}\subseteq Mor(X,Y)_\beta$ and $\exists\sigma_1,...,\sigma_p\in J_\beta=\prod\limits_{x\in X}H_{\beta(x)}$ such that $\sigma_if_0=f_i\ (1\leq i\leq p)$. Define $\sigma_i=(h_i(x))_{x\in X}$. Then $(\sigma_i f_0)(x)=(\sigma_i)_xf_0(x)=h_i(x)f_0(x)=f_i(x)$ for any $x\in X$ and hence $\sigma_if_0=f_i$. Since $f_0\in Mor(X,Y)_\beta$ it follows that $f_i\in Mor(X,Y)_\beta\ (1\leq i\leq p)$ and hence $\{f_0,...,f_p\}\subseteq Mor(X,Y)_\beta$. \hfill$\Box$
}
\Proof{\textbf{of Lemma 2.28} \hspace{0,1cm}Let $\{x_0,...x_q\}$ be a local frame in X. Since $f_0 \in Mor(X,Y)_\beta$ it follows that $\{f_0(x_0),...,f_0(x_q)\}$ is a b-frame in Y such that $b\geq\beta(x_0),...,\beta(x_q)$. Hence $\exists\sigma_1,...,\sigma_q\in H_b$ such that $\sigma_j f_0(x_0)=f_0(x_j)\ (1\leq j\leq q)$. Since $\{f_0,...,f_p\}$ is a $\beta$-frame in $Mor(X,Y)$ $\exists\tau_1,...,\tau_p\in J_\beta$ such that $\tau_i f_0=f_i\ (1\leq i\leq p)$. Let $i\in \{0,...,p\}, j\in\{0,...,q\}$. Then (notice that $f_0(x_j)\in Y_{\beta(x_j)}\cap Y_b$) \begin{align*}&\hspace{0.45cm}(H_{\beta(x_j)\leq b}((\tau_i)_{x_j})\sigma_j)f_0(x_0)\\&=H_{\beta(x_j)\leq b}((\tau_i)_{x_j})(\sigma_jf_0(x_0))\\&=H_{\beta(x_j)\leq b}((\tau_i)_{x_j})f_0(x_j)\\&=(\tau_i)_{x_j}f_0(x_j)\\&=f_i(x_j),\end{align*} i.e $H_b$ acts transitively on $\{f_i(x_j)|0\leq i\leq p, 0\leq j\leq q\}$ and hence $\{f_i(x_j)|0\leq i\leq p, 0\leq j\leq q\}$ is a b-frame.\hfill $\Box$
}
\Proof{\textbf{of Lemma 2.29} \hspace{0,1cm}Let Y be an infimum action, X be a global action and $g:X\rightarrow Y$ a morphism of global actions. We have to show that g is Z-normal for any global actions Z, i.e. if Z is a global action and $\{f_0,...,f_p\}$ is a local frame in $Mor(Z,X)$ then $\{gf_0,...,gf_p\}$ is a local frame in $Mor(X,Y)$. Let Z be a global action and $\{f_0,...,f_p\}$ be a local frame in $Mor(Z,X)$. From Lemma 2.27 it follows that $\{f_0(z),...,f_p(z)\}$ is a local frame in X $\forall z\in Z$. Since $g:X\rightarrow Y$ is a morphism of global actions it follows that $\{(gf_0)(z),...,(gf_p)(z)\}$ is a local frame in Y $\forall z\in Z$. For any $z\in Z$ the set $\{b\in\Psi|\{(gf_0)(z),...,(gf_p)(z)\}\textnormal{ is a b-frame}\}$ has an initial element $b_z$, since Z is an infimum action. Define \begin{align*}\beta:Z&\rightarrow\Psi\\ z&\mapsto b_z\ .\end{align*} We show now that $gf_0\in Mor(Z,Y)_\beta$, i.e. \begin{enumerate}[(a)]
\item $(gf_0)(z)\in Y_{\beta(z)}=Y_{b_z}\ \forall z\in Z$, and
\item if $\{z_0,...,z_q\}$ is a local frame in z then $\{(gf_0)(z_0),...,(gf_0)(z_q)\}$ is a c-frame for some $c\in\Psi$ such that $c\geq\beta(z_i)=b_{z_i}\ (0\leq i\leq q)$. 
\end{enumerate}Clearly $\{(gf_0)(z),...,(gf_p)(z)\}$ is a $b_z$-frame for any $z\in Z$. Thus $(gf_0)(z)\in Y_{b_z}\ \forall z\in Z$. We show that (b) holds. Let $\{z_0,...z_q\}$ be a local frame in Z. By Lemma 2.28 $\{f_i(z_j)|0\leq i\leq p, 0\leq j\leq q\}$ is a local frame in X and hence $\{(gf_i)(z_j)|0\leq i\leq p, 0\leq j\leq q\}$ is a c-frame for some $c\in\Psi$. This implies $\{(gf_0)(z_0),...,(gf_0)(z_q)\}$ is a c-frame and $\{(gf_0)(z_j),...,(gf_p)(z_j)\}$ is a c-frame $\forall j\in\{0,...,q\}$, whence $c\geq b_{z_j}\ (0\leq j\leq q)$. Thus (b) holds and hence $gf_0\in Mor(Z,Y)_\beta$. Since $\{(gf_0)(z),...,(gf_p)(z)\}$ is a $\beta(z)=b_z$-frame $\forall z\in Z$, it follows from Lemma 2.27 that $\{gf_0,...,gf_p\}$ is a $\beta$-frame.\hfill$\Box$
}
\Proof{\textbf{of Lemma 2.30} \hspace{0,1cm}First we show that if Z is an infimum action then for any global action X, $Mor(X,Z)$ is an infimum action. Let Z be an infimum action, i.e. if $U\subseteq Z$ is a finite and nonempty subset then the set $\{c\in\Theta|U\textnormal{ is a c-frame}\}$ is either empty or contains an initial element. Let X be a global action. We must show that $Mor(X,Z)$ is an infimum action, i.e. if $V\subseteq Mor(X,Z)$ is a finite and nonempty subset then the set $\{\gamma\in\Theta_{(X,Z)}|V\textnormal{ is a }\gamma\textnormal{-frame}\}$ is either empty or contains an initial element. Let $V=\{f_0,...,f_p\}\subseteq Mor(X,Z)$ be a finite and nonemtpy subset. Suppose $\{\gamma\in\Theta_{(X,Z)}|V\textnormal{ is a }\gamma\textnormal{-frame}\}$ is nonempty. Let $\delta\in\{\gamma\in\Theta_{(X,Z)}|V\textnormal{ is a }\gamma\textnormal{-frame}\}$. Then V is a $\delta$-frame and hence it follows from Lemma 2.27 that $\{f_0(x),...f_p(x)\}$ is a $\delta(x)$-frame $\forall x\in X$. Since Z is an infimum action,  $\{c\in\Theta|\{f_0(x),...f_p(x)\}\textnormal{ is a c-frame}\}$ contains an initial element $c_x$ for any $x\in X$. Define \begin{align*}\delta':X&\rightarrow\Theta\\ x&\mapsto c_x\ .\end{align*}. Clearly $\delta(x)\geq c_x=\delta'(x)\ \forall x\in X\Leftrightarrow\delta\geq\delta'$. We will show now that $f_0\in Mor(X,Z)_{\delta'},i.e.$ \begin{enumerate}[(a)]
\item $f_0(x)\in Z_{\delta'(x)}\ \forall x\in X$, and
\item if $\{x_0,...,x_q\}$ is a local frame in X then $\{f_0(x_0),...,f_0(x_q)\}$ is a d-frame for some $d\in\Theta$ such that $d\geq\delta'(x_i)\ (0\leq i\leq q)$.
\end{enumerate}Clearly $\{f_0(x),...,f_p(x)\}$ is a $c_x=\delta'(x)$-frame $\forall x\in X$. Thus $f_0(x)\in Z_{\delta'(x)} \forall x\in X$. (b) holds because of the following. Let $\{x_0,...,x_q\}$ be a local frame in X. Since $\{f_0,...,f_p\}$ is a $\delta$-frame in $Mor(X,Z)$, it follows from Lemma 2.28 that $\{f_i(x_j)|0\leq i\leq p, 0\leq j\leq q\}$ is a d-frame for some $d\in\Theta$. This implies that $\{f_0(x_0),...,f_0(x_q)\}$ is a d-frame and that $\{f_0(x_j),...,f_p(x_j)\}$ is a local frame $\forall j\in\{0,...,q\}$, whence $d\geq c_{x_j}=\delta'(x_j)\ (0\leq j \leq q)$. Thus (b) holds and hence $f_0\in Mor(X,Z)_{\delta'}$. Since $\{f_0(x),...,f_p(x)\}$ is a $c_x=\delta'(x)$-frame $\forall x\in X$, it follows from Lemma 2.27 that $\{f_0,...,f_p\}$ is a $\delta'$-frame. Thus $\delta'\in\{\gamma\in\Theta_{(X,Z)}|V\textnormal{ is a }\gamma\textnormal{-frame}\}$.\\

We show now that if Z is a strong infimum action and the relation on $\Theta$ is transitive then for any global action X, $Mor(X,Z)$ is a strong infimum action and the relation on $\Theta_{(X,Z)}$ is transitive. Let Z be a strong infimum action, i.e. if $\Delta\subseteq\Theta$ denotes a finite subset and $U\subset Z$ denotes a finite and nonempty subset such that $U\cap Z_\gamma\neq\emptyset\ \forall\gamma\in\Delta$ then the set $\{c\in\Theta_{\geq\Delta}|U\textnormal{ is a c-frame}\}$ is either empty or contains an initial element.  Let the relation on $\Theta$ be transitive and let X be a global action. We must show that $Mor(X,Z)$ is a strong infimum action, i.e. if $\Gamma\subseteq\Theta_{(X,Z)}$ is a finite subset and $V\subseteq Mor(X,Z)$ is a finite and nonempty subset such that $V\cap Mor(X,Z)_\delta\neq\emptyset\ \forall\delta\in\Gamma$ then the set $\{d\in(\Theta_{(X,Z)})_{\geq\Gamma}|V\textnormal{ is a }d\textnormal{-frame}\}$ is either empty or contains an initial element, and that the relation on $Mor(X,Z)$ is transitive. Let $\alpha,\beta,\gamma\in \Theta_{(X,Z)}$ such that $\alpha\leq\beta$, i.e. $\alpha(x)\leq\beta(x)\ \forall x\in X$, and $\beta\leq\gamma$, i.e. $\beta(x)\leq\gamma(x)\ \forall x\in X$. Since the relation on $\Theta$ is transitive it follows that $\alpha(x)\leq\gamma(x)\ \forall x\in X$ which is equivalent to $\alpha\leq\gamma$. Thus the relation on $\Theta_{(X,Z)}$ is transitive. Let $\Gamma\subseteq\Theta_{(X,Z)}$ be a finite subset and $V=\{f_0,...,f_p\}\subseteq Mor(X,Z)$ be a finite and nonemtpy subset such that $V\cap Mor(X,Z)_\delta\neq\emptyset\ \forall\delta\in\Gamma$. Suppose $\{d\in(\Theta_{(X,Z)})_{\geq\Gamma}|V\textnormal{ is a }d\textnormal{-frame}\}$ is nonempty. Let $e\in\{d\in(\Theta_{(X,Z)})_{\geq\Gamma}|V\textnormal{ is a }d\textnormal{-frame}\}$. Then V is an $e$-frame and hence it follows from Lemma 2.27 that $\{f_0(x),...,f_p(x)\}$ is an $e(x)$-frame $\forall x\in X$. Let $\Delta(x)=\{\delta(x)|\delta\in\Gamma\}$. Let $x\in X$ and $\epsilon\in\Delta(x)$. Then $\epsilon=\delta(x)$ for some $\delta\in\Gamma$. Since $V\cap Mor(X,Z)_{\delta'}\neq\emptyset\ \forall\delta'\in\Gamma$, $\exists i\in\{0,...,p\}$ such that $f_i\in Mor(X,Z)_\delta$. It follows that $f_i(x)\in Z_{\delta(x)}=Z_\epsilon$. Hence if $x \in X$ then $\{f_0(x),...,f_p(x)\}\cap Z_\epsilon\neq\emptyset\ \forall\epsilon\in\Delta(x)$. Since Z is a strong infimum action,  $\{c\in\Theta_{\geq\Delta(x)}|\{f_0(x),...,f_p(x)\}\textnormal{ is a c-frame}\}$ contains an initial element $c_x$ for any $x\in X$. Define \begin{align*}e':X&\rightarrow\Theta\\ x&\mapsto c_x\ .\end{align*}Since $e\geq\delta\ \forall\delta\in\Gamma$ it follows that $e(x)\geq\delta(x)\ \forall\delta\in\Gamma\ \forall x\in X$. Hence $e(x)\in\Theta_{\geq\Delta(x)}\ \forall x\in X$. Since $\{f_0(x),...,f_p(x)\}$ is an $e(x)$-frame for any $x\in X$, it follows $e(x)\in\{c\in\Theta_{\geq\Delta(x)}|\{f_0(x),...,f_p(x)\}\textnormal{ is a c-frame}\}\ \forall x\in X$. Thus $e(x)\geq c_x=e'(x)\ \forall x\in X\Leftrightarrow e\geq e'$. We will show now that $f_0\in Mor(X,Z)_{e'},i.e.$ \begin{enumerate}[(a)]
\item $f_0(x)\in Z_{e'(x)}\ \forall x\in X$, and
\item if $\{x_0,...,x_q\}$ is a local frame in X then $\{f_0(x_0),...,f_0(x_q)\}$ is an f-frame for some $f\in\Theta$ such that $f\geq e'(x_i)\ (0\leq i\leq q)$.
\end{enumerate}Clearly $\{f_0(x),...,f_p(x)\}$ is a $c_x=e'(x)$-frame $\forall x\in X$. Thus $f_0(x)\in Z_{e'(x)} \forall x\in X$. (b) holds because of the following. Let $\{x_0,...,x_q\}$ be a local frame in X. Since $\{f_0,...,f_p\}$ is an e-frame in $Mor(X,Z)$, it follows from Lemma 2.28 that $\{f_i(x_j)|0\leq i\leq p, 0\leq j\leq q\}$ is an f-frame for some $f\in\Theta$ such that $f\geq e(x_0),...,e(x_p)$. This implies that $\{f_0(x_0),...,f_0(x_q)\}$ is an f-frame and that $\{f_0(x_j),...,f_p(x_j)\}$ is a local frame $\forall j\in\{0,...,q\}$. Since $e'(x_i)\geq\delta(x_i)\ \forall\delta\in\Gamma\ \forall i\in\{0,...,q\}$ it follows, by the transitivity of the relation on $\Theta$, that $f\geq \delta(x_i)\ \forall\delta\in\Gamma\ \forall i\in\{0,...,q\}$. Thus $f\in\Theta_{\geq\Delta(x_i)}\ (0\leq i\leq q)$. It follows that $f\in\{c\in\Theta_{\geq\Delta(x)}|\{f_0(x_i),...,f_p(x_i)\}\textnormal{ is a c-frame}\}\ $ $(0\leq i\leq q)$ and hence $f\geq c_{x_j}=\gamma'(x_j)\ (0\leq j \leq q)$. Thus (b) holds and hence $f_0\in Mor(X,Z)_{e'}$. Since $\{f_0(x),...,f_p(x)\}$ is a $c_x=e'(x)$-frame $\forall x\in X$, it follows from Lemma 2.27 that $\{f_0,...,f_p\}$ is a $e'$-frame. Thus $e'\in\{d\in\Theta_{(X,Z)}|V\textnormal{ is a }d\textnormal{-frame}\}$. It remains to show that $e'\geq\delta\ \forall\delta\in\Gamma$. Let $\delta\in\Gamma$. Since $e'(x)=c_x\in\{c\in\Theta_{\geq\Delta(x)}|\{f_0(x),...,f_p(x)\}\textnormal{ is a c-frame}\}$ for any $x\in X$, it follows that  $e'(x)\geq\delta(x)\ \forall x\in X \Leftrightarrow e'\geq\delta$. Thus $e'\in\{d\in(\Theta_{(X,Z)})_{\geq\Gamma}|V\textnormal{ is a }d\textnormal{-frame}\}$.\hfill$\Box$

\section{Homotopy}
\Definition{Let X and Y denote global actions and let $f, g:X\rightarrow Y$ denote morphisms. Let L denote the line action defined in 2.3. We say that f is \textit{homotopic} to g if there is a morphism $H:X\times L\rightarrow Y$ and elements $n_-, n_+\in\mathbb{Z}, n_-\leq n_+$ such that $H\overset{\circ}{\iota}_n=f\ \forall n\leq n_-$ and $H\overset{\circ}{\iota}_n=g\ \forall n\geq n_+$ where \begin{align*}\overset{\circ}{\iota}_n:X&\rightarrow X\times L\\ x&\mapsto (x,n).\end{align*}
}
\Definition[homotopy]{Let X, Y denote global actions. A morphism $H:X\times L\rightarrow Y$ is called a \textit{homotopy}\index{homotopy} if $\exists \textnormal{ elements } n_-\leq n_+\in\mathbb{Z}$ such that $H\overset{\circ}{\iota}_n=H\overset{\circ}{\iota}_{n_+}\ \forall n\geq n_+$ and $H\overset{\circ}{\iota}_n=H\overset{\circ}{\iota}_{n_-}\ \forall n\leq n_-$.}
\Definition{Let X denote a global action. A morphism $f:L\rightarrow X$ \textit{stabilizes on the left}, if $\exists n_-\in\mathbb{Z}$ such that $f(n)=f(n_-)\ \forall n\leq n_-$. f \textit{stabilizes on the right}, if $\exists n_+\in\mathbb{Z}$ such that $f(n)=f(n_+)\ \forall n\geq n_+$. A pair $(n_+,n_-)$ as above is called a \textit{stabilisation pair} for f\index{morphism!stabilisation pair of a -}.} 
\Definition[path]{Let X denote a global action. A morphism $f:L\rightarrow X$ which stabilizes on the left and on the right is called a \textit{path}\index{path} in X.}
\Definition{Let X denote a global action and $\omega:L\rightarrow X$ a path. If $\omega$ is not constant then there is a smallest integer $lus(\omega)$ called the \textit{least upper stabilization}\index{least upper stabilization} of $\omega$ such that $\omega(n)=\omega(lus(\omega))\ \forall n\geq lus(\omega)$ and a largest integer $gls(\omega)$ called the \textit{greatest lower stabilization}\index{greatest lower stabilization} of $\omega$ such that $\omega(n)=\omega(gls(\omega))\ \forall n\leq gls(\omega)$. Suppose $\omega$ is constant then we define $lus(\omega)=0$ and $gls(\omega)=0$. $\omega(lus(\omega))$ is called the \textit{terminal point}\index{terminal point} of $\omega$ ($\textit{term}(\omega)$) and $\omega(gls(\omega))$ is called the \textit{initial point}\index{initial point} of $\omega$ ($\textit{init}(\omega)$).}
\Definition[loop]{A \textit{loop}\index{loop} in X is a path $\omega$ in X whose initial point and terminal point are equal, i.e. if $n_+=lus(\omega)$ and $n_+=gls(\omega)$ then $\omega(n_-)=\omega(n_+)$. If $x=\textit{init}(\omega)$ then $\omega$ is called a \textit{loop at $x$}. The set of all loops at x is denoted $\Omega(x)$.
\Definition{Let $\omega$ and $\omega'$ denote loops in a global action X. We say that $\omega$ is \textit{loop homotopic} to $\omega'$ \begin{enumerate}[(a)]
\item if $\omega$ is homotopic to $\omega'$ in the usual sense, i.e. there is a morphism $H:L\times L\rightarrow X$ and integers $n_-, n_+\in\mathbb{Z}, n_-\leq n_+$, such that $H\overset{\circ}{\iota}_n=\omega\ \forall n\leq n_-$ and $H\overset{\circ}{\iota}_n=\omega'\ \forall n\geq n_+$, and
\item $H\overset{\circ}{\iota}_n$ is a loop for any $n\in\mathbb{Z}$.
\end{enumerate}
}
\Definition[composition of paths\index{path!composition of -s}]{Let $\omega$ and $\omega'$ denote paths in a global action X such that $\textit{term}(\omega)=\textit{init}(\omega')$. If $\omega$ is not constant then their \textit{composition} $\omega'\omega$ is defined as follows. \[\omega'\omega(n)=\begin{cases}
\omega(n),  & \text{if }n\leq lus(\omega),\\
\omega'(n-lus(\omega)+gls(\omega')), & \text{if }n\geq lus(\omega).
\end{cases}\]
The composition is well defined, since $\omega(lus(\omega))=\textit{term}(\omega)=\textit{init}(\omega')=\omega'(gls(\omega'))$ $=\omega'(lus(\omega)-lus(\omega)+gls(\omega'))$. Clearly, if $\omega'$ is constant then $\omega'\omega=\omega$. If $\omega$ is constant then define $\omega'\omega=\omega'$.
}
\Remarkc{Clearly loops in $\Omega(x)$ compose and their composition is a loop at $x$. The definition of composition in 3.8 is concocted so that the composition law induced on $\Omega(x)$ is a monoid with identity element the constant function $x$.
}
\Definition[inverse of a path\index{path!inverse of a -}]{Let $\omega$ be a path in a global action X. The map $\omega^{-1}:l \rightarrow X$ defined by \[\omega^{-1}(n)=\begin{cases}
\textit{term}(\omega),  & \text{if }n<gls(\omega),\\
\omega(lus(\omega)-n+gls(\omega)), & \text{if }gls(\omega)\leq n\leq lus(\omega),\\
\textit{init}(\omega), &\text{if }n>lus(\omega),
\end{cases}\] is called the \textit{inverse} of $\omega$. It is clearly a path.}
\Remark{Let $\omega$ and $\omega'$ be paths in a global action $X$ such that $\textit{term}(\omega)=\textit{init}(\omega')$. Then $\omega\omega'$ and $\omega^{-1}$ also are paths in X.}
\Definition[path connected]{Let $\{G_\alpha\curvearrowright X_\alpha|\alpha\in\Phi\}$ denote a global action on X. Two points $x,y\in X$ are called \textit{path connected} if there is a path $\omega$ in X such that $\textit{init}(\omega)=x$ and $\textit{term}(\omega)=y$. If x and y are path connected we write $x\sim y$. X is called \textit{(path) connected}\index{global action!connected -} if any two points $x,y\in X$ are path connected.}
\DefinitionLemma{Let X be a global action. Then $\sim$ is an equivalence relation on $\bigcup\limits_{\alpha\in\Phi}X_\alpha$. By definition $\pi_0(X)$ is the set of all equivalence classes of $\sim$ on $\bigcup\limits_{\alpha\in\Phi}X_\alpha$. An equivalence class of $\sim$ is called a \textnormal{path component}\index{path!- component} of X.}
\Proof{We must check that $\sim$ is reflexive, symmetric and transitive.\\ 
\\
$\bullet$ $\sim$ is reflexive.\\
Let $x\in X$. Let $\omega$ be the constant path in X at x. Then $\textit{init}(\omega)=\textit{term}(\omega)=x$ and hence $x\sim x$.\\
\\
$\bullet$ $\sim$ is symmetric.\\
Let $x,y\in X$ such that $x\sim y$. Then there is a path $\omega$ in X such that $\textit{init}(\omega)=x$ and $\textit{term}(\omega)=y$. Since $\textit{init}(\omega^{-1})=\textit{term}(\omega)=y$ and $\textit{term}(\omega^{-1})=\textit{init}(\omega)=x$, it follows that $y\sim x$.\\
\\
$\bullet$ $\sim$ is transitive.\\
Let $x,y,z\in X$ such that $x\sim y$, $y\sim z$. Then there are paths $\omega, \omega'$ in X such that $\textit{init}(\omega)=x$, $\textit{term}(\omega)=y$, $\textit{init}(\omega')=y$ and $\textit{term}(\omega)=z$. Clearly $\omega'\omega$ is a path in X such that $\textit{init}(\omega'\omega)=x$ and $\textit{term}(\omega'\omega)=z$. Hence $x\sim z$. $\Box$ 
}
\Definition{Let X,Y denote global actions, $H:X\times L\rightarrow Y$ denote a homotopy. Suppose $\exists m,n\in\mathbb{Z}$, $m\neq n$, such that $H\overset{\circ}{\iota}_m\neq H\overset{\circ}{\iota}_n$. Then there is a largest $n_-(H)\in\mathbb{Z}$ such that $H\overset{\circ}{\iota}_n=H\overset{\circ}{\iota}_{n_-(H)}\ \forall n\leq n_-(H)$ and a smallest $n_+(H)\in\mathbb{Z}$ such that $H\overset{\circ}{\iota}_n=H\overset{\circ}{\iota}_{n_+(H)}\ \forall n\geq n_+(H)$. If $H\overset{\circ}{\iota}_n=H\overset{\circ}{\iota}_m\ \forall m,n\in\mathbb{Z}$ then we define $n_-(H)=n_+(H)=0$. $n_-(H)$ is called the \textit{greatest lower stabilization} of H (\textit{$gls(H)$}) and $n_+(H)$ is called the \textit{least upper stabilization} of H (\textit{$lus(H)$}). $H\overset{\circ}{\iota}_{n_-(H)}$ is called the \textit{initial value} of H ($\textit{init}(H)$) and $H\overset{\circ}{\iota}_{n_+(H)}$ is called the \textit{terminal value} of H ($\textit{term}(H)$).}
\Definition[composition of homotopies\index{homotopy!composition of -ies}]{Let $H_1,H_2:X\times L\rightarrow Y$ be homotopies such that $\textit{term}(H_1)=\textit{init}(H_2)$. Then their \textit{composition} $H_2H_1$ is defined as follows.\[H_2H_1((x,n))=\begin{cases}
H_1((x,n)),  & \text{if }n\leq lus(H_1),\\
H_2((x,n-lus(H_1)+gls(H_2))), & \text{if }n\geq lus(H_1).
\end{cases}\]}
\Definition[inverse of a homotopy\index{homotopy!inverse of a -}]{Let $H:X\times L\rightarrow Y$ denote a homotopy. The map $H^{-1}:X\times L\rightarrow Y$ defined by \[H^{-1}(x,n)=\begin{cases}
H(x,gls(H)),  & \text{if }n<gls(H),\\
H(x,lus(H)-n+gls(H)), & \text{if }gls(\omega)\leq n\leq lus(\omega),\\
H(x,lus(H)), &\text{if }n>lus(H),
\end{cases}\] is called the \textit{inverse} of $H$. It is clearly a homotopy.}
}
\Definition[base point\index{base point}]{Let X denote a global action. Giving X a \textit{base point} means fixing some $*\in X$. If $*$ is the base point of X we write sometimes $X_*$ instead of X.}
\Definition[end point stable homotopy\index{homotopy!end point stable -}]{Let X denote a global action with base point $*$ and  $H:L\times L\rightarrow X$ a homotopy. If for any $n\in Z$, $H\overset{\circ}{\iota}_n$ is a loop whose initial point(=terminal point) is $*$, then H is called an \textit{end point stable homotopy}.}
\DefinitionLemma{Let X denote a global action with base point $*$. Define a relation $\sim$ on $\Omega(*)$ by $\omega\sim\omega'\Leftrightarrow$ $\omega$ and $\omega'$ are end point stable homotopic, i.e. there is a homotopy $H:L\times L\rightarrow X$ such that $\forall n\in \mathbb{Z}\ H\overset{\circ}{\iota}_n$ is a loop at $*$ and for all n sufficiently small (resp. large) $H\overset{\circ}{\iota}_n=\omega$ (resp. $H\overset{\circ}{\iota}_n =\omega'$). Then $\sim$ is an equivalence relation on $\Omega(*)$ and, by definition, $\pi_1(X_*)$ \index{fundamental group} is the set of all equivalence classes of $\sim$. If X is pathwise connected and $|\pi_1(X_*)|=1$ , X is called \textnormal{simply connected}\index{global action!simply connected -}. }
\Proof{We have to check that $\sim$ is reflexive, symmetric and transitive.\\ 
\\
$\bullet$ $\sim$ is reflexive.\\
Let $\omega\in \Omega(*)$. Define a map $H:L\times L\rightarrow X$ by $h(m,n)=\omega(m)\ \forall m,n\in\mathbb{Z}$. One checks easily that H is an end point stable homotopy such that $\textit{init}(H)=\textit{term}(H)=\omega$. Hence $\omega\sim\omega$.\\
\\
$\bullet$ $\sim$ is symmetric.\\
Let $\omega,\omega'\in \Omega(*)$ such that $\omega\sim\omega'$. Then there is an end point stable homotopy $H:L\times L\rightarrow X$ such that $\textit{init}(H)=\omega$ and $\textit{term}(H)=\omega'$. Since $\textit{init}(H^{-1})=\textit{term}(H)=\omega'$ and $\textit{term}(H^{-1})=\textit{init}(H)=\omega$, it follows that $\omega'\sim\omega$.\\
\\
$\bullet$ $\sim$ is transitive.\\
Let $\omega,\omega', \omega''\in X$ such that $\omega\sim\omega'$, $\omega'\sim\omega''$. Then there are end point stable homotopies  $H, H':L\times L\rightarrow X$  such that $\textit{init}(H)=\omega$, $\textit{term}(H)=\omega'$, $\textit{init}(H')=\omega'$ and $\textit{term}(H')=\omega''$. Clearly $H'H:L\times L\rightarrow X$ is an end point stable homotopy such that $\textit{init}(H'H)=\omega$ and $\textit{term}(H'H)=\omega''$. Hence $\omega\sim\omega''$. $\Box$ 
}
\Lemma{Let X denote a global action with basepoint $*$. If $\omega\in\Omega(x)$, let $[\omega]$ denote the equivalence class of loops $\omega'$ which are end point stable homotopic to $\omega$. Define \begin{align*}\circ:\pi_1(X_*)\times\pi_1(X_*)&\rightarrow\pi_1(X_*)\\ ([\omega],[\tau])&\mapsto[\tau\omega]\ .\end{align*} Then $\circ$ is well defined and $(\pi_1(X_*), \circ)$ is a group.}
\Proof{First we show that $\circ$ is well defined. Let $A,B\in\pi_1(X_*)$, $\omega,\omega'\in A$ and $\tau,\tau'\in B$. We must show that $[\omega\tau]=[\omega'\tau']$, i.e. $\omega\tau\sim\omega'\tau'$. Since $\omega\sim\omega'$ and $\tau\sim\tau'$, there are end point stable homotopies $H_1,H_2:L\times L\rightarrow X$ such that $\textit{init}(H_1)=\omega$, $\textit{term}(H_1)=\omega'$, $\textit{init}(H_2)=\tau$ and $\textit{term}(H_2)=\tau'$. Define a map $H:L\times L\rightarrow X$ by $H\overset{\circ}{\iota}_n=H_1\overset{\circ}{\iota}_nH_2\overset{\circ}{\iota}_n\ (n\in\mathbb{Z})$. Obviously H is an end point stable homotopy such that $\textit{init}(H)=\omega\tau$ and $\textit{term}(H)=\omega'\tau'$ and hence $\omega\tau\sim\omega'\tau'$. We show now that $(\pi_1(X_*), \circ)$ is a group. Since $\Omega(*)$ is a monoid by Remark 3.9, it follows that $(\pi_1(X_*), \circ)$ is a monoid. Hence we only have to show the existence of inverse elements. Let $A=[\omega]\in\pi_1(X_*)$. Define a function $f_k:L\rightarrow X$ by \[f_k(n)=\begin{cases}
\omega^{-1}\omega(n),  & \text{if }n<lus(\omega)-k,\\
\omega^{-1}\omega(lus(\omega)-k), & \text{if }lus(\omega)-k\leq n\leq lus(\omega)+k,\\
\omega^{-1}\omega(n), & \text{if }n>lus(\omega)+k.
\end{cases}\] ($0\leq k\leq lus(\omega)-gls(\omega)$). One checks easily that $f_{k-1}$ is end point stable homotopic to $f_{k}$ by a homotopy $H_k\ (0\leq k\leq lus(\omega)-gls(\omega)-1)$. Clearly $f_0=\omega^{-1}\omega$ is end point stable homotopic to $f_{lus(\omega)-gls(\omega)}=\textnormal{constant loop at }*$ by $H_{lus(\omega)-gls(\omega)-1}...H_1H_0$. Hence $[\omega]\circ[\omega^{-1}]=[\omega^{-1}\omega]=[\tau]$. $\Box$}

\section{Coverings}
\Definition{Let $p:Y\rightarrow X$ denote a morphism of global actions. p is called a \textit{covering morphism}\index{morphism!covering -} if and only if for any local frame $\{x_0,...,x_n\}$ in X and any $y_0\in p^{-1}(x_0)$ there is a unique local frame $\{y_0,y_1...,y_n\}$ in Y such that $p(y_i)=x_i\ \forall i\in \{1,...,n\}$.}
\Definition{Let X denote a global action and $x\in X$. Let $star(x)=\bigcup\limits_{x\in X_\alpha}G_\alpha x$. We give $star(x)$ the structure of a global action as follows. Let $\Phi_{star(x)}=\{\alpha\in\Phi|\ x\in X_\alpha\}$, $(G_{star(x)})_\alpha=G_\alpha$ and $star(x)_\alpha=G_\alpha x$. We let $(G_{star(x)})_\alpha$ act on $star(x)_\alpha$ in the obvious way. Note: Let $x,x_1,...,x_p$ be a finite set of elements in X. Then $x,x_1,...,x_p$ is a local frame in X $\Leftrightarrow x,x_1,...,x_p\in star(x)$ and is a local frame there also.
\Lemma{Let X,Y be global actions and let $p:Y\rightarrow X$ be a map. Then p is a covering morphism if and only if $y\in Y$, $p|_{star(y)}:star(y)\rightarrow star(p(y))$ is an isomorphism of global actions.}
\Proof{Suppose $p:Y\rightarrow X$ is a covering morphism. We have to check several things.\\
\\
$\bullet$ $p|_{star(y)}(y')\in star(p(y))\ \forall y'\in star(y)$.\\
Let $y'\in star(y)$. Since $y'\in star(y)$, there is a $\beta\in\Psi_{star(y)}$ and an $h\in H_\beta$ such that $hy=y'$. Obviously $\{y,y'\}$ is a local frame in Y. Since p is a morphism of global actions, it follows that $\{p(y),p(y')\}$ is a local frame in X. Thus $p(y')\in star(x)$.\\
\\
$\bullet$ $p|_{star(y)}$ is a morphism of global actions.\\
Let $\{y_0,...,y_q\}$ be a local frame in $star(y)$. It follows that $\{y_0,...,y_q\}$ is a local frame in Y. Since p is a morphism of global actions $\{p|_{star(y)}(y_0),...,p|_{star(y)}(y_q)\}=\{p(y_0),...,p(y_q)\}$ is a local frame in X. Hence $\{p|_{star(y)}(y_0),...,p|_{star(y)}(y_q)\}$ is a local frame in star(p(y)).\\
\\
$\bullet$ $p|_{star(y)}$ is injective.\\
Let $y_1, y_2\in star(y)$ such that $p(y_1)=p(y_2)$. Since $\{p(y), p(y_1)\}$ is a local frame in X, there is a unique local frame $\{y, y'\}$ in Y such that $p(y')=p(y_1)$. Since $\{y, y_1\}$ and $\{y, y_2\}$ are both liftings of the same local frame in X, it follows that they are equal and hence $y_1=y_2$.\\
\\
$\bullet$ $p|_{star(y)}$ is surjective.\\
Let $x\in star(p(y))$, i.e. $\{p(y), x\}$ is a local frame in X. Since p is a covering morphism, there is a local frame $\{y, y'\}$ in Y such that $p(y')=x$. Since $\{y, y'\}$ is a local frame, $y'\in star(y)$.\\
\\
$\bullet$ Since $p|_{star(y)}$ is bijective, it suffices to show that $(p|_{star(y)})^{-1}:star(p(y))\rightarrow star(y)$ is a morphism of global actions.\\
Let $\{x_0,...,x_q\}$ be a local frame in $star(p(y))$. Then $\{x_0,...,x_q\}\in G_\alpha p(y)$ for some $\alpha\in\Phi_{star(p(y))}$. Thus $\{p(y),x_0,...,x_q\}$ is a local frame in X. Hence there is a (unique) local frame $\{y,y_0,...,y_q\}$ such that $p(y_i)=x_i\ (0\leq i\leq q)$. Obviously $y_0,...,y_q\in star(y)$. Thus $(p|{star(y)})^{-1}$ is a morphism.\\
\\
Suppose now $p|_{star(y)}:star(y)\rightarrow star(p(y))$ is an isomorphism of global actions $\forall y\in Y$. We show first that p is a morphism. Let $\{y_0,...,y_q\}$ be a local frame in Y. Then $\{y_0,...,y_q\}$ is a local frame in $star(y_0)$. Since $p|_{star(y_0)}:star(y_0)\rightarrow star(p(y_0))$ is a morphism, $\{p(y_0),...,p(y_q)\}$ is a local frame in star $(p(y_0))$ and thus also in X. We show now that p satisfies the unique local frame lifting property. Let $\{x_0,...x_q\}$ be a local frame in X and let $y_0\in p^{-1}(x_0)$. We must show that there is a unique local frame $\{y_0,y_1,...,y_q\}$ in Y such that $p(y_i)=x_i\ \forall i\in \{1,...,n\}$. \\
\\
$\bullet$ existence\\
Obviously $\{x_0,...x_q\}\subseteq star(x_0)$. Since $p|_{star(y)}:star(y_0)\rightarrow star(x_0)$ is an isomorphism of global actions, $(p|_{star(y_0)})^{-1}:star(x_0)\rightarrow star(y_0)$ is an isomorphism of global actions. Hence $\{(p|_{star(y_0)})^{-1}(x_0),(p|_{star(y_0)})^{-1}(x_1),...,(p|_{star(y_0)})^{-1}(x_q)\}=\{y_0,(p|_{star(y_0)})^{-1}(x_1),...,(p|_{star(y_0)})^{-1}(x_q)\}$ is a local frame in Y.\\
\\
$\bullet$ uniqueness\\
Suppose $\{y_0,y_1,...,y_q\}$ is a local frame in Y such that $p(y_i)=x_i\ \forall i\in \{1,...,n\}$. Then $\{y_0,y_1,...,y_q\}\subseteq star(y_0)$. Since $p|_{star(y_0)}:star(y_0)\rightarrow star(x_0)$ is injective, it follows that $y_i=(p|_{star(y_0)})^{-1}(x_i)\ (1\leq i\leq q)$. $\Box$ 
}
\Lemma[Unique Path Lifting Property (UPLP)\index{Unique Path Lifting Property}]{Let $p:Y\rightarrow X$ be a covering morphism of global actions, $\omega$ a path in X and $(n_-,n_+)$ be a stabilisation pair for $\omega$. Let $y\in p^{-1}(\textit{init}(\omega))$. Then there is a unique path $\tilde\omega$ in Y such that $\textit{init}(\tilde\omega)=y$. $(n_-,n_+)$ is also a stabilisation pair for $\tilde\omega$.
}
\Proof{See \cite{atlas}, p.154, proof of Lemma 10.4. $\Box$}
\Definition[fixed end point homotopy\index{homotopy!fixed end point -}]{Let X denote a global action, $x, y\in X$ and $\omega, \omega'$ be paths in X such that $\textit{init}(\omega)=\textit{init}(\omega')=a$ and $\textit{term}(\omega)=\textit{term}(\omega')=b$. A morphism $H:L\times L\rightarrow X$ such that $\exists n_-, n_+\in \mathbb{Z}, n_-\leq n_+$ such that $\overset{\circ}{\iota}_n=\omega\ \forall n\leq n_-$, $\overset{\circ}{\iota}_n=\omega'\ \forall n\geq n_+$, $\overset{\circ}{\iota}_n(m)=a\ \forall m\leq n_-\forall n\in\mathbb{Z}$ and $\overset{\circ}{\iota}_n(m)=b\ \forall m\geq n_+\forall n\in\mathbb{Z}$. }
\Lemma{\index{Homotopy Lifting Theorem}Let $p:Y\rightarrow X$ be a covering of global actions. Let $\omega, \omega'$ be paths in X. Let $H:L\times L\rightarrow X$ be a fixed end point homotopy from $\omega$ to $\omega'$. Let $y\in p^{-1}(\textit{init}(\omega))$. Then the lifts $\tilde\omega$ and $\tilde\omega'$ of $\omega$ resp. $\omega'$ such that $\textit{init}(\tilde\omega)=\textit{init}(\tilde\omega')=y$ are homotopic by a fixed end point homotopy $\tilde H$ (this implies $\textit{term}(\tilde\omega)=\textit{term}(\tilde\omega')$). 
}
\Proof{See \cite{atlas}, p.154 f., proof of Lemma 10.5. $\Box$}
\Theorem[Lifting Criterion\index{Lifting Criterion}]{Let $(X,x_0), (Y,y_0)$ and $(Z,z_0)$ denote global actions with base points. Let $p:Y\rightarrow X$ denote a covering morphism such that $p(y_0)=x_0$ and $f:Z\rightarrow X$ a morphism such that $f(z_0)=x_0$. Let $p_*:\pi_1(Y,y_0)\rightarrow\pi_1(X,x_0)$ and $f_*:\pi_1(Z,z_0)\rightarrow\pi_1(X,x_0)$ be the induced maps. Suppose that Z is path connected. Then f lifts to a morphism $\tilde f:Z\rightarrow Y$ such that $\tilde f(z_0)=y_0$ if and only if $f_*(\pi_1(Z,z_0))\subseteq p_*(\pi_1(Y,y_0))$.
}
\Proof{
$\Rightarrow$: Suppose f lifts to a morphism $\tilde f:Z\rightarrow Y$ such that $\tilde f(z_0)=y_0$, i.e. there is a morphism $\tilde f:Z\rightarrow Y$ such that $\tilde f(z_0)=y_0$ and $p\tilde f=f$. Then \begin{align*} f_*(\pi_1(Z,z_0))=(p\tilde f)_*(\pi_1(Z,z_0))=p_*(\tilde f_*(\pi_1(Z,z_0)))\subseteq p_*(\pi_1(Y,y_0)).\end{align*}
$\Leftarrow$: Suppose $f_*(\pi_1(Z,z_0))\subseteq p_*(\pi_1(Y,y_0))$. Define a map $\tilde f:Z\rightarrow Y$ as follows. Let $z\in Z$. Since Z is connected there is a path $\omega$ in Z such that $\textit{init}(\omega)=z_0$ and $\textit{term}(\omega)=z$. Clearly $y_0\in p^{-1}(f\omega(gls(\omega)))=p^{-1}(f(\textit{init}(\omega)))=p^{-1}(f(z_0))=p^{-1}(x_0)$. By the UPLP $f\omega$ lifts uniquely to a path $\tilde \omega$ in Y such that $\textit{init}(\tilde \omega)=y_0$. Define $\tilde f(z)=\textit{term}(\tilde \omega)$. There are several things to check.\\
\\
$\bullet$ $\tilde f$ is well defined.\\
Let $\omega$ and $\omega'$ be paths in Z such that $\textit{init}(\omega)=\textit{init}(\omega')=z_0$ and $\textit{term}(\omega)=\textit{term}(\omega')=z$. Then $(f\omega)^{-1}(f\omega')$ is a loop in X at $x_0$, since $\textit{init}((f\omega)^{-1}(f\omega'))=\textit{init}(f\omega')=f(\textit{init}(\omega'))=f(z_0)=x_0$ and $\textit{term}((f\omega)^{-1}(f\omega'))=\textit{term}(f\omega)^{-1}=\textit{init}(f\omega)=f(\textit{init}(\omega))=f(z_0)=x_0$. Clearly $[(f\omega)^{-1}(f\omega')]\in f_*(\pi_1(Z,z_0))\subseteq p_*(\pi_1(Y,y_0))$. Hence there is a loop $\tau$ in Y at $y_0$ such that $p_*([\tau])=[(f\omega)^{-1}(f\omega')]\Leftrightarrow [p\tau]=[(f\omega)^{-1}(f\omega')]$, i.e. $p\tau$ and $(f\omega)^{-1}(f\omega')$ are homotopic by a homotopy $H:L\times L\rightarrow X$. By the UPLP $(f\omega)^{-1}(f\omega')$ lifts uniquely to a path $\lambda$ in Y with the stabilisation pair $(gls(\omega'), lus(\omega')+lus(\omega)-gls(\omega))$ such that $\textit{init}(\lambda)=y_0$. By Lemma 4.6 $\textit{term}(\lambda)=\textit{term}(\tau)=y_0$. Set $n_0=lus(\omega')+gls(\omega)-lus(\omega)$. For all $n\in\mathbb{Z}$ define \[\mu(n)=\begin{cases}
\lambda(n),  & \text{if }n\leq lus(\omega'),\\
\lambda(lus(\omega')), & \text{if }n\geq lus(\omega')
\end{cases}\] and \[\nu(n)=\begin{cases}
y_0,  & \text{if }n\leq gls(\omega'),\\
\lambda(lus(\omega')+lus(\omega)-gls(\omega)-n+lus(\omega')), & \text{if }gls(\omega')\leq n \leq n_0,\\
\lambda(lus(\omega')), & \text{if }n\geq n_0.
\end{cases}\] Obviously $\mu=\tilde\omega'$ and $\nu=\tilde\omega$. Hence $\textit{term}(\tilde\omega)=\textit{term}(\nu)=\lambda(lus(\omega'))=\textit{term}(\mu)=\textit{term}(\tilde\omega')$.\\
\\
$\bullet$ $p\tilde f=f$. \\
Let $z\in Z$. Let $\omega$ be a path in Z such that $\textit{init}(\omega)=z_0$ and $\textit{term}(\omega)=z$. Let $\tilde\omega$ be the unique lift of $f\omega$ such that $\textit{init}\tilde\omega=y_0$. It follows that \begin{align*}(p\tilde f)(z)=p(\tilde f(z))=p(\textit{term}(\tilde\omega))=p(\tilde\omega(lus(\omega))=p(f\omega(lus(\omega))=f(z).\end{align*}\\
\\
$\bullet$ $\tilde f$ is a morphism of global actions.\\
Let $\{z_1,...,z_q\}$ be a local frame in Z. We have to show that $\{\tilde f(z_1),...,\tilde f(z_q)\}$ is a local frame in Y. Let $\omega_1$ be an arbitrary path in Z such that  $\textit{init}(\omega_1)=z_0$ and $\textit{term}(\omega_1)=z_1$. For any $i\in\{2,...,q\}$, $n\in\mathbb{Z}$ define  \[\omega_i(n)=\begin{cases}
\omega_1(n),  & \text{if }n\leq lus(\omega_1),\\
z_i, & \text{if }n>lus(\omega_1).
\end{cases}\] Then \[\tilde\omega_i(n)=\begin{cases}
\tilde\omega_1(n),  & \text{if }n\leq lus(\omega_1),\\
\tilde f(z_i), & \text{if }n>lus(\omega_1), 
\end{cases}\] for any $i\in\{2,...,q\}$, $n\in\mathbb{Z}$. Hence $\tilde f(z_2),...,\tilde f(z_q)\in star(\tilde f(z_1))$. Since f is a morphism of global actions $\{f(z_1),...,$ $f(z_q)\}$ a local frame in X. Since p is a covering morphism there is a local frame $\{\tilde f(z_1),y_2,...,y_q\}$ in Y such that $p(y_i)=f(z_i)\ (2\leq i\leq q)$. But since $p|_{star(\tilde f(z_1))}$ is bijectice by Lemma 4.3, $y_i=\tilde f(z_i)\ (2\leq i\leq q)$.\\ 
\\
$\bullet$ $\tilde f(z_0)=y_0$.\\
Let $\omega$ be the constant path at $z_0$. Let $\tilde \omega$ be the unique lifting of $\omega$ such that $\textit{init}(\tilde\omega)=y_0$. By definition of $\tilde f$, $\tilde f(z_0)=\textit{term}(\tilde\omega)$. Since $(gls(\omega), lus(\omega))=(0,0)$ is a stabilisatioin pair for $\tilde \omega$, $\tilde\omega$ is the constant path at $y_0$ and hence $\textit{term}(\tilde\omega)=y_0$. $\Box$  
}
\Theorem{Let $(X,x_0)$ denote a connected global action with base point and let H be a subgroup of $\pi_1(X,x_0)$. Then there is a connected global action $X_H$, an element $\tilde x_0$ and a covering morphism $p_H:X_H\rightarrow X$ such that ${p_H}_*(\pi_1(X_H,\tilde x_0))=H$. 
}
\Proof{
Let $\{X_\alpha\curvearrowright G_\alpha|\alpha\in\Phi\}$ denote a connected global action with base point $x_0$ and let H be a subgroup of $\pi_1(X,x_0)$. Let $\Gamma X$ denote the set of all paths in X starting at $x_0$. Define a relation $\sim_H$ on $\Gamma X$ by $\omega\sim_H\omega'\ \Leftrightarrow\ \textit{term}(\omega)=\textit{term}(\omega')$ and $[(\omega')^{-1}\omega]\in H$. One checks easily that $\sim_H$ is an equivalence relation. Let $X_H$ denote the set of all equivalence classes of $\sim_H$. Define a function $p_H: X_H\rightarrow X$ by $p_h([\omega])=\textit{term}(\omega)$. We will give now $A_H$ the structure of a global action. Define $\Phi_H=\Phi$, $(X_H)_\alpha=\{a\in X_H|p_H(a)\in(X_H)_\alpha\}$ and $(G_H)_\alpha=G_\alpha$. Define an action of $(G_H)_\alpha$ on $(X_H)_\alpha$ by $\sigma [\omega]$=$[\sigma\omega]$ where \[(\sigma\omega)(n)=\begin{cases}
\omega(n),  & \text{if }n\leq lus(\omega),\\
\sigma\textit{term}(\omega), & \text{if }n>lus(\omega).
\end{cases}\]We will show that the global action $X_H$ defined above is connected and that $p_H:X_H\rightarrow X$ is a covering morphism such that ${p_H}_*(\pi_1(X_H,\tilde x_0))=H$, where $\tilde x_0$ is the equivalence class of the constant path at $x_0$.\\
\\
$\bullet$ $p_H:X_H\rightarrow X$ is a morphism.\\
Let $\{[\omega_0],...,[\omega_p]\}$ be a local frame in $X_H$, i.e. $\{[\omega_0],...,[\omega_p]\}\subseteq (X_H)_\alpha$ for some $\alpha\in\Phi_H=\Phi$ and $\exists \sigma_1,...,\sigma_p\in (G_H)_\alpha=G_\alpha$ such that $\sigma_i[\omega_0]=[\omega_i]\ (1\leq i \leq p)$. We have to show that $\{p_H([\omega_0]),...,p_H([\omega_p])\}=\{\textit{term}(\omega_0),...,\textit{term}(\omega_p)\}$ is a local frame in X. Clearly $\{\textit{term}(\omega_0),...,\textit{term}(\omega_p)\}\subseteq X_\alpha$, since $\{[\omega_0],...,[\omega_p]\}\subseteq (X_H)_\alpha$. Let $i\in\{1,...,p\}$. Since $[\sigma_i\omega_0]=\sigma_i[\omega_0]=[\omega_i]$, it follows that $\sigma_i\textit{term}(\omega_0)=\textit{term}(\sigma_i\omega_0)=\textit{term}(\omega_i)$. Hence $\{p_H([\omega_0]),...,p_H([\omega_p])\}=\{\textit{term}(\omega_0),...,\textit{term}(\omega_p)\}$ is a local frame in X.\\
\\
$\bullet$ $p_H:X_H\rightarrow X$ is a covering morphism.\\
Let $\{x_0,...,x_p\}$ be a local frame in X and $a_0=[\omega]\in p_H^{-1}(x_0)$, i.e. $\{x_0,...,x_p\}\subseteq X_\alpha$ for some $\alpha\in\Phi$ and $\exists g_1,...,g_p\in G_\alpha$ such that $g_i x_0=x_i\ (0\leq i\leq p)$ . We must show that there is a unique local frame $\{a_0,a_1,...,a_p\}$ in $X_H$ such that $p_H(a_i)=x_i\ (1\leq i\leq p)$. First we show existence. Set $a_i=g_ia_0$. Then obviously $\{a_0,a_1,...,a_p\}$ is a local frame in $X_H$ with the property above. Assume now that $\{a_0,b_1,...,b_p\}$ is a local frame in $X_H$ with the property above. Then $\{a_0,b_1,...,b_p\}\subseteq (X_H)_\beta$ for some $\beta\in\Phi_H=\Phi$ and $\exists h_1,...,h_p\in(G_H)_\beta$ such that $h_ia_0=b_i\ (1\leq i\leq p)$. Since $h_i\textit{term}(\omega)=p_H([h_i\omega])=p_H(h_i[\omega])=p_H(h_ia_0)=p_H(b_i)=x_i\ (1\leq i\leq p)$ and analogously $g_i\textit{term}(\omega)=x_i$, it follows that $h_i\omega=g_i\omega\ (1\leq i \leq p)$. Hence $b_i=h_ia_0=h_i[\omega]=[h_i\omega]=[g_i\omega]=g_i[\omega]=g_ia_0=a_i\ (1\leq i\leq p)$.
}\\
\\
$\bullet$ $X_H$ is connected.\\
Let $[\omega]\in X_H$. By the UPLP $\omega$ lifts (uniquely) to a path $\tau$ in $X_H$ starting at $\tilde x_0$. Clearly, if $n\in \{1,...,lus(\omega)\}$ , then $\tau(n)=[\rho_n]$, where \[\rho_n(m)=\begin{cases}
\omega(m),  & \text{if }m\leq n,\\
\omega(n), & \text{if }m>n.
\end{cases}\] Since $\textit{term}(\rho_{lus(\omega)})=\textit{term}(\omega)$ and $[(\rho_{lus(\omega)})^{-1}\omega]=[\omega^{-1}\omega]=[\textit{constant path at }x_0]\in H$, it follows that $\rho_{lus(\omega)}\sim_H\omega$ and hence $\textit{term}(\tau)=[\rho_{lus(\omega)}]=[\omega]$. Thus $\tilde x_0$ and $[\omega]$ are path connected and since $[\omega]$ was an arbitrary element of $X_H$, $X_H$ is connected.\\
\\
$\bullet$ ${p_H}_*(\pi_1(X_H,\tilde x_0))=H$.\\
Let $[l]\in\pi_1(X_H,\tilde x_0)$. Clearly $p_H l$ is a loop in $X$ at $x_0$ and $l$ is its unique lift starting at $\tilde x_0$. It follows that $[c]=\tilde x_0=\textit{term}(l)=[p_H l]$ (see above), where $c$ is the constant path in $X$ at $x_0$. Hence $[(p_H l)c^{-1}]\in H$. Since $\textit{term}((p_H l)c^{-1})=\textit{term}(p_H l)$ and $[((p_h l)c^{-1})^{-1}(p_h l)]=[(p_H l)^{-1}(p_H l)]=[c]\in H$, it follows that ${p_H}_*([l])=[p_H l]=[(p_h l)c^{-1}]\in H$. $\Box$
}
\Corollary {A connected global action has a simply connected covering space.}
\Proof{Let $(X,x_0)$ denote a connected global action. First we show that ${p_H}_*:\pi_1(X_H,\tilde x_0)\rightarrow H$ is injective. Let $[l],[l']\in \pi_1(X_H,\tilde x_0)$ such that ${p_H}_*([l])={p_H}_*([l'])$. By Lemma 4.6, $[l]=[l']$ and hence ${p_H}_*$ is injective. Let H={[c]}, where $c$ is the constant path in $X$ at $x_0$. It follows that $|\pi_1(X_H,\tilde x_0)|\leq1$, i.e. $X_H$ is simply connected. $\Box$  
}
\newpage
\addcontentsline{toc}{section}{References} 
\nocite*{} 
\bibliographystyle{plain} 
\bibliography{literatur_diplom} 
\addcontentsline{toc}{section}{Index}
\index{covering|see{morphism}}
\index{product of global actions|see{global action}}
\index{morphism space|see{global action}}
\index{strong infimum condition|see{global action}}
\index{stabilisation pair|see{morphism}}
\printindex
\end{document}

%% file: Meta.tex
\newcommand{\titel}{Global Actions}

\newcommand{\art}{Diploma Thesis															}

\newcommand{\autor}{Raimund Preusser}
\newcommand{\studienbereich}{mathematics}
\newcommand{\matrikelnr}{1633126}
\newcommand{\betreuer}{Prof. Anthony Bak}

\newcommand{\jahr}{2010}



%% file: Deckblatt.tex
\thispagestyle{plain}
\begin{titlepage}

\begin{center}
\vspace*{1cm}
\huge{\textbf{\textsc{\titel}}}\\[1.5ex]
\LARGE{\textbf{\art}}\\[1.5ex]
\vspace{4cm}
\includegraphics[scale=0.25]{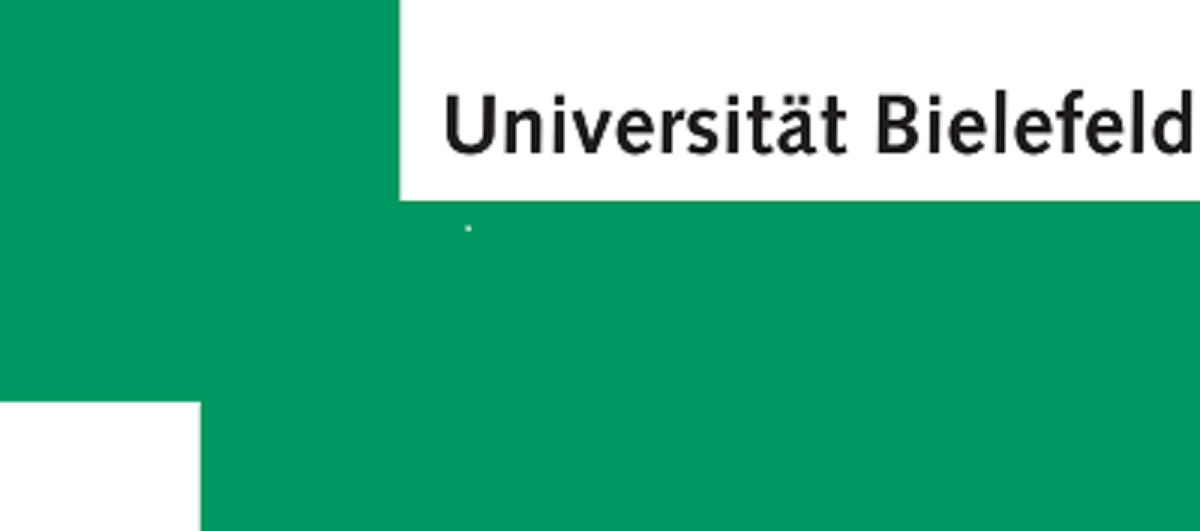}\\[3ex]
\vspace{4cm}
\normalsize
\begin{tabular}{lc}\\
 submittet by:	 & \quad \autor\\[1.2ex]
 field of study: & \quad \studienbereich\\[1.2ex]
 matr. number: & \quad \matrikelnr\\[1.2ex]
 thesis adviser:         & \quad \betreuer\\[1.2ex]
\end{tabular}\\
\vspace*{2.3cm}
\hrule
\vspace{0,2cm}
\hspace{\fill}\copyright\ \jahr\\[1.5ex]

\end{center}


\end{titlepage}